\documentclass[11pt]{article}
\usepackage{amssymb,amsmath,amsthm, amsfonts}

\def\sqr#1#2{{\vcenter{\vbox{\hrule height.#2pt
              \hbox{\vrule width.#2pt height#1pt \kern#1pt \vrule width.#2pt}
          \hrule height.#2pt}}}}
\def\signed #1{{\unskip\nobreak\hfil\penalty50
          \hskip2em\hbox{}\nobreak\hfil#1
          \parfillskip=0pt \finalhyphendemerits=0 \par}}
\def\endpf{\signed {$\sqr69$}}
\def\sqr#1#2{{\vcenter{\vbox{\hrule height.#2pt
              \hbox{\vrule width.#2pt height#1pt \kern#1pt \vrule width.#2pt}
              \hrule height.#2pt}}}}
\def\signed #1{{\unskip\nobreak\hfil\penalty50
              \hskip2em\hbox{}\nobreak\hfil#1
              \parfillskip=0pt \finalhyphendemerits=0 \par}}
\def\endpf{\signed {$\sqr69$}}
\def\3n{\negthinspace \negthinspace \negthinspace }
\def\2n{\negthinspace \negthinspace }
\def\1n{\negthinspace }

\def\={\buildrel \triangle \over =}

\def\O{\Omega}

\def\q{\quad}

\def\max{\mathop{\rm max}}

\def\exp{\mathop{\rm exp}}
\def\sup{\mathop{\rm sup}}

\def\|{\Big |}
\def\({\Big (}
\def\){\Big )}
\def\[{\Big[}
\def\]{\Big]}
\def\be{\begin{equation}}
\def\bel{\begin{equation}\label}
\def\ee{\end{equation}}
\def\bt{\begin{theorem}}
\def\bcd{\begin{condition}}
\def\ecd{\end{condition}}
\def\et{\end{theorem}}
\def\bc{\begin{corollary}}
\def\ec{\end{corollary}}
\def\bde{\begin{definition}}
\def\ede{\end{definition}}
\def\bl{\begin{lemma}}
\def\el{\end{lemma}}
\def\bp{\begin{proposition}}
\def\ep{\end{proposition}}
\def\br{\begin{remark}}
\def\er{\end{remark}}
\def\ba{\begin{array}}
\def\ea{\end{array}}
\def\ed{\end{document}}

\def\square#1{\vbox{\hrule\hbox{\vrule height#1%
     \kern#1\vrule}\hrule}}
\def\rectangle#1#2{\vbox{\hrule\hbox{\vrule height#1%
     \kern#2\vrule}\hrule}}

\font\tenbb=msbm10 \font\sevenbb=msbm7 \font\fivebb=msbm5

\newfam\bbfam
\scriptscriptfont\bbfam=\fivebb \textfont\bbfam=\tenbb
\scriptfont\bbfam=\sevenbb

\newtheorem{lemma}{Lemma}[section]
\newtheorem{remark}{Remark}[section]

\newtheorem{theorem}{Theorem}[section]
\newtheorem{corollary}{Corollary}[section]

\newtheorem{definition}{Definition}[section]
\newtheorem{proposition}{Proposition}[section]
\newtheorem{condition}{Condition}[section]

\makeatletter
   
   \@addtoreset{equation}{section}
\makeatother

\begin{document}

\title{ Mean-Field Backward Stochastic Differential Equations and
Related Partial Differential Equations \bf\footnote{The work has
been supported by the NSF of P.R.China (No. 10701050; 10671112),
Postdoctoral Science Foundation of Shanghai grant (No. 06R214121)
and National Basic Research Program of China (973 Program) (No.
2007CB814904; 2007CB814906).}}

\author{ Rainer Buckdahn\\
{\small D\'{e}partement de Math\'{e}matiques, Universit\'{e} de
Bretagne Occidentale,}\\
 {\small 6, avenue Victor-le-Gorgeu, B.P. 809, 29285 Brest
cedex, France.}\\
{\small{\it E-mail: Rainer.Buckdahn@univ-brest.fr.}}\\
 Juan Li\\
{\small Department of Mathematics, Shandong University at Weihai, Weihai 264200, P. R. China.;}\\
{\small Institute of Mathematics, School of Mathematical Sciences,
Fudan University, Shanghai 200433.}\\
{\small {\it E-mail: juanli@sdu.edu.cn.}}\\
Shige Peng\\
{\small School of Mathematics and System Sciences, Shandong University, Jinan 250100, P. R. China.}\\
{\small {\it E-mail: peng@sdu.edu.cn.}}  }

\date{}

\maketitle \noindent{\bf Abstract}\hskip4mm

In~\cite{BLP} the authors obtained Mean-Field backward stochastic
differential equations (BSDE) associated with a Mean-field
stochastic differential equation (SDE) in a natural way as limit of
some highly dimensional system of forward and backward SDEs,
corresponding to a large number of ``particles'' (or ``agents'').
The objective of the present paper is to deepen the investigation of
such Mean-Field BSDEs by studying them in a more general framework,
with general driver, and to discuss comparison results for them. In
a second step we are interested in partial differential equations
(PDE) whose solutions can be stochastically interpreted in terms of
Mean-Field BSDEs. For this we study a Mean-Field BSDE in a Markovian
framework, associated with a Mean-Field forward equation. By
combining classical BSDE methods, in particular that of ``backward
semigroups" introduced by Peng~\cite{Pe1}, with specific arguments
for Mean-Field BSDEs we prove that this Mean-Field BSDE describes
the viscosity solution of a nonlocal PDE. The uniqueness of this
viscosity solution is obtained for the space of continuous functions
with polynomial growth. With the help of an example it is shown that
for the nonlocal PDEs associated to Mean-Field BSDEs one cannot
expect to have uniqueness in a larger space of continuous functions.


\newpage
\section{\large{Introduction}}

\hskip1cm Classical mathematical mean-field approaches can be met in
different fields, so in Statistical Mechanics and Physics (for
instance, the derivation of Boltzmann or Vlasov equations in the
kinetic gas theory) and in Quantum Mechanics and Quantum Chemistry
(for instance, the density functional models or also Hartree and
Hartree-Fock type models). In a recent series of papers
(see~\cite{LL} and their papers cited therein) Lasry and Lions
extended the field of such mean-field approaches to problems in
Economics and Finance, and also to the theory of stochastic
differential games. Recently Buckdahn, Li and Peng~\cite{BLP}
studied a special mean field problem in a purely stochastic approach
and deduced a new kind of BSDEs which they called Mean-Field BSDEs.
In order to be more precise, we consider the stochastic dynamics
$X=(X_t)_{t\in[0,T]}$ of a particle \be
\begin{array}{rcl}
dX_t&=&\sigma(t,X_t,X^1_t,\dots,X^N_t)dB_t+b(t,X_t,X^1_t,\dots,
X^N_t)dt,\, t\in[0,T],\\
X_0&=&x\in R^n, \end{array}\ee which is governed by a
$d$-dimensional Brownian motion $B=(B_t)_{ t\in[0,T]}$ and
influenced by the positions $X^{1}, \dots,X^{N}$ of $N$ other
particles that are supposed to be mutually independent of each
other, also independent of the driving Brownian motion $B$, and of
the same law as $X$. The time $T>0$ is an arbitrarily fixed horizon.
In~\cite{BLP} it is shown that, under appropriate assumptions, SDE
(1.1) has a unique solution.

For the study of the limit behavior $N\rightarrow +\infty$ of the
above SDE of rank $N$, the authors of~\cite{BLP} assumed that the
coefficients $\sigma=\sigma^N$\ and $b=b^N$\ of SDE (1.1) of rank
$N$ have the form
$$\displaystyle\sigma^N(t,x,x^1,\dots,x^N)=\frac{1}{N}\sum_{i=1}^N
\sigma(t,x^{i},x),\ \ b^N(t,x,x^1,\dots,x^N)=\frac{1}{N}\sum_{i=1}^N
b(t,x^{i},x),$$ \noindent $(t,x,(x^1,\dots,x^N))\in [0,T]\times
R^n\times R^{n\times N}$. The functions $\sigma:[0,T]\times
R^n\times R^n\rightarrow R^{n\times d}$, $b:[0,T]\times R^n\times
R^n\rightarrow R^n$\ are supposed to be bounded, measurable  and,
moreover, Lipschitz in the space variables, uniformly with respect
to the time variable. Denoting by
$(X^{(N)},X^{(N),1},\dots,X^{(N),N}):=(X,X^1,\dots,X^N)$  the
solution of SDE (1.1) with the coefficients $\sigma^N$ and $b^N$, it
is shown in~\cite{BLP} that the the processes $X^{(N)}$ converge
uniformly on the time interval $[0,T]$, in $L^2$, to the unique
continuous $\mathbb{F}^B$-adapted solution $\hat{X}$ of the limit
equation (Mean-Field SDE)\be
d\hat{X}_t=E\left[\sigma(t,\hat{X}_t,\mu)\right]_{\mu=\hat{X}_t}dB_t
+E\left[\sigma(t,\hat{X}_t,\mu) \right]_{\mu=\hat{X}_t}dt,\,
t\in[0,T],\, \hat{X}_0=x\ee ($\mathbb{F}^B$ is the filtration
generated by $B$).

With the SDE of rank $N$ the authors of \cite{BLP} associated a BSDE
of the form \be\begin{array}{rcl}
dY_t&=&-f^N(t,X_t,X^1_t,\dots,X^N_t,Y_t,Y^1_t,\dots,Y^N_t,Z_t,Z^1_t,
\dots,Z^N_t)dt+Z_tdB_t,\ t\in[0,T],\\
Y_T&=&\Phi^N(X_T,X^1_T,\dots,X^N_T),\\
\end{array} \ee in which, in the same spirit as for the forward SDE of rank
$N$, $(X^{i},Y^{i},Z^{i}),1\le i\le N,$ is a sample of $N$ mutually
independent triplets which are independent of the driving Browning
motion $B$ and obey the same law as the solution $(X,Y,Z)$; recall
that $(X,X^1,\dots,X^N)$ is the solution of the above SDE (1.1) with
coefficients $\sigma^N$ and $b^N$. The coefficients $f^N,\Phi^N$ for
the above BSDE are introduced in the same spirit as $\sigma^N,b^N$:

$\displaystyle
f^N(t,x,x^1,\dots,x^N,y,y^1,\dots,y^N,z,z^1,\dots,z^N)
=\frac{1}{N}\sum_{i=1}^Nf(t,x^{i},x,y^{i},y,z^{i},z),$

$\displaystyle\Phi^N(x,x^1,\dots,x^N)=\frac{1}{N}\sum_{i=1}^N
\Phi(x^{i},x)$,

\smallskip

\noindent $t\in[0,T],\  x,x^j\in R^n,y,y^j\in R,z,z^j\in R^d,\, 1\le
j\le N$, where $f:[0,T]\times R^n\times R^n\times R\times R\times
R^d\times R^d\rightarrow R$ and $\Phi:R^n\times R^n\rightarrow R$
are bounded, measurable and Lipschitz in the variables
$(x',x,y',y,z',z)$, uniformly with respect to $t\in [0,T]$.
In~\cite{BLP} it is proved that, with
$$(Y^{(N)},Y^{(N),1},\dots,Y^{(N), N},Z^{(N)},Z^{(N),1},\dots, Z^{(N),
N}):=(Y,Y^1,\dots,Y^N,Z,Z^1,\dots, Z^N)$$ denoting the solution of
the above BSDE (1.3), the couple of processes $(Y^{(N)},Z^{(N)})$,
interpreted as random variable with values in $C([0,T];R)\times
L^2([0,T];R^d),$ converges to the unique square integrable,
$\mathbb{F}^B$-progressively measurable solution of the limit
equation (called Mean-Field BSDE) \be\begin{array}{rcl}\displaystyle
d\hat{Y}_t&=&-E\left[f(t,\hat{X}_t,\mu,\hat{Y}_t,\lambda,\hat{Z}_t,\zeta)
\right]_{\mu=\hat{X}_t,\lambda=\hat{Y}_t,\zeta=\hat{Z}_t}dt+\hat{Z}_tdB_t,\,\ t\in[0,T],\\
Y_T&=&E\left[\Phi(\hat{X}_T,\mu)\right]_{\mu=\hat{X}_T},\end{array}\ee
where $\hat{X}$\ is the solution of Mean-Field SDE (1.2).

The objective of the present paper is to deepen the investigation of
the above Mean-Field BSDE. In a first leg we study the existence and
uniqueness for Mean-Field BSDEs in a rather general setting, with
drivers which, in difference to~\cite{BLP}, are not necessarily
deterministic coefficients. In addition to the existence and the
uniqueness also the comparison principle for this new type of BSDEs
is discussed and illustrated by examples.

The main objective of the paper concerns the study of Mean-Field
problems in a Markovian setting. To be more precise, we investigate
Mean-Field BSDEs associated with Mean-Field forward SDEs and partial
differential equations (PDEs) whose solutions are described by them.
The system dynamics we investigate is given by the following SDE \be
  \left\{
  \begin{array}{rcl}
  dX_s^{t,\zeta}&=&E[b(s,X_s^{0,x_0}, \mu)]_{\mu=X_s^{t,\zeta}}ds+
  E[\sigma(s,X_s^{0,x_0},\mu)]_{\mu=X_s^{t,\zeta}}dB_s,\ s\in[t,T],\\
  X_t^{t,\zeta}&=&\zeta.
  \end{array}
  \right.
  \ee
Precise assumptions on the coefficients $b:[0,T]\times
\mathbb{R}^n\times \mathbb{R}^n\rightarrow \mathbb{R}^n$ and
$\sigma:[0,T]\times \mathbb{R}^n\times \mathbb{R}^n\rightarrow
\mathbb{R}^{n\times d}$ are given in the following sections.

With the above SDE we associate the BSDE: \be \left\{
\begin{array}{rcl}
-dY_s^{t,\zeta}&=&E[f(s,X_s^{0,x_0},\mu,
Y_s^{0,x_0},\lambda,\zeta)]_{\mu=X_s^{t,\zeta},\lambda=Y_s^{t,\zeta},
\zeta=Z_s^{t,\zeta}}ds
-Z_s^{t,\zeta}dB_s,\ s\in [t, T],\\
Y_T^{t,\zeta}&=&E[\Phi(X_T^{0,x_0},\mu)]_{\mu=X_T^{t,\zeta}}.\\
\end{array}
\right. \ee Under the assumptions on $f$ and $\Phi$ that are
introduced in Section 5, the above BSDE has a unique solution
$(Y_s^{t,x},Z_s^{t,x})_{s\in[t,T]}$ and we can define the
deterministic function
\begin{equation}
u(t,x)=Y_t^{t,x}.
\end{equation}
We prove that $u(t,x)$\ is the unique viscosity solution in $C_p([0,
T]\times {\mathbb{R}}^n)$ to the following nonlocal PDE \be\left
\{\begin{array}{ll}
 &\!\!\!\!\! \frac{\partial }{\partial t} u(t,x) + Au(t,x)
+ E[f(t, X_t^{0, x_0}, x, u(t, X_t^{0, x_0}), u(t,x),
Du(t,x).E[\sigma(t, X_t^{0, x_0},
 x)])]=0,\\
 & \hskip 4.3cm (t,x)\in [0,T)\times {\mathbb{R}}^n ,  \\
 &\!\!\!\!\!  u(T,x) =E[\Phi (X_T^{0, x_0}, x)], \hskip0.5cm   x \in
 {\mathbb{R}}^n,
 \end{array}\right.
\ee with
$$ Au(t,x):=\frac{1}{2}tr(E[\sigma(t, X_t^{0, x_0}, x)]
 E[\sigma(t, X_t^{0, x_0},x)]^{T}D^2u(t,x))+Du(t,x).E[b(t,
 X_t^{0, x_0},x)].$$
In particular, it is shown that the space $C_p([0, T]\times
{\mathbb{R}}^n)$ is the optimal space in which the uniqueness can be
got.

Our paper is organized as follows. Section 2 recalls some elements
of the theory of BSDEs which are needed in what follows. Section 3
investigates the properties of general Mean-Field BSDEs. We first
prove the uniqueness and existence of the solution of Mean-Field
BSDE (Theorem 3.1) but also the comparison theorem (Theorem 3.2) and
the converse comparison theorem (Theorem 3.3). Similarly we
investigate Mean-Field Forward SDEs in Section 4. In Section 5 we
investigate decoupled Mean-Field Forward-Backward SDEs (FBSDEs).
Their value function $u$\ (see (5.4)) turns out to be a
deterministic function which is Lipschitz in $x$\ (see (5.5)) and
$\frac{1}{2}$-H\"{o}lder continuous in $t$\ (Theorem 5.2). Moreover,
it satisfies the dynamic programming principle (DPP) (see 5.10). A
key element in the proof of the DPP is the use of Peng's backward
semigroups (see~\cite{Pe1}). We change slightly its definition; this
allows to shorten the argument for the proof that $u$\ is a
viscosity solution of the associated PDE (Theorem 6.1). Finally, the
uniqueness of the viscosity solution in the space of continuous
functions with polynomial growth is proved in Section 7.

\section{ {\large Preliminaries}}

  \hskip1cm Let $\{B_t\}_{t\geq0}$\ be a d-dimensional standard Brownian
  motion defined over some complete probability space $(\Omega,{\cal{F}},
  P).$\ By ${\mathbb{F}}=\{{\mathcal{F}}_s,\ 0\leq s \leq
T\}$\ we denote the natural filtration generated by $\{B_s\}_{0\leq
s\leq T}$\ and augmented by all P-null sets, i.e.,
$${\mathcal{F}}_s=\sigma\{B_r, r\leq s\}\vee {\mathcal{N}}_P,\ s\in [0, T], $$
where $ {\cal{N}}_P$ is the set of all P-null subsets and $T > 0$\ a
fixed real time horizon.  For any
    $n\geq 1,$\ $|z|$ denotes the Euclidean norm of $z\in
    {\mathbb{R}}^{n}$. We also shall introduce the following both spaces of
    processes which are used frequently in what follows:
\vskip0.2cm
    ${\cal{S}}_{\mathbb{F}}^2(0, T; {\mathbb{R}}):=\{(\psi_t)_{0\leq t\leq T}
    \mbox{ real-valued}\ {\mathbb{F}} \mbox{-adapted c\`{a}dl\`{a}g
    process}:\\ \mbox{ }\hskip6cm
    E[\sup\limits_{0\leq t\leq T}| \psi_{t} |^2]< +\infty \}; $
    \vskip0.2cm

   ${\cal{H}}_{\mathbb{F}}^{2}(0,T;{\mathbb{R}}^{n}):=\{(\psi_t)_{0\leq t\leq T}
   \ {\mathbb{R}}^{n}\mbox{-valued}\ {\mathbb{F}} \mbox{-progressively
   measurable process}:\\ \mbox{ }\hskip6cm
     \parallel\psi\parallel^2_{2}=E[\int^T_0| \psi_t| ^2dt]<+\infty \}. $
\vskip0.2cm
 Let us now consider a function $g:\Omega\times[0,T]\times {\mathbb{R}} \times {\mathbb{R}}^{d}
\rightarrow {\mathbb{R}} $ with the property that $(g(t, y,
z))_{t\in [0, T]}$\ is ${\mathbb{F}}$-progressively measurable for
each $(y,z)$\ in ${\mathbb{R}} \times {\mathbb{R}}^{d}$, and which
is assumed to satisfy the following standard assumptions throughout
the paper:
 \vskip0.2cm

(A1) There exists a constant $C\ge 0$  such that, \mbox{dtdP-a.e.,
for all}\ $y_{1}, y_{2}\in {\mathbb{R}},\ z_{1}, z_{2}\in
{\mathbb{R}}^d,\\ \mbox{ }\hskip4cm   |g(t, y_{1}, z_{1}) - g(t,
y_{2}, z_{2})|\leq C(|y_{1}-y_{2}| + |z_{1}-z_{2}|);$
 \vskip0.2cm

(A2) $g(\cdot,0,0)\in {\cal{H}}_{\mathbb{F}}^{2}(0,T;{\mathbb{R}})$.
\vskip0.2cm

 The following result on BSDEs is by now well known, for its proof the reader is referred to
 Pardoux and Peng~\cite{PaPe}.
 \bl Under the assumptions (A1) and (A2), for any random variable
 $\xi\in L^2(\O, {\cal{F}}_T,$ $P),$ the
BSDE
 \be y_t = \xi + \int_t^Tg(s,y_s,z_s)ds - \int^T_tz_s\,
dB_s,\q 0\le t\le T, \label{BSDE} \ee
 has a unique adapted solution
$$(y_t, z_t)_{t\in [0, T]}\in {\cal{S}}_{\mathbb{F}}^2(0, T;
{\mathbb{R}})\times
{\cal{H}}_{\mathbb{F}}^{2}(0,T;{\mathbb{R}}^{d}). $$ \el
 In what follows, we always assume that the driving coefficient $g$\ of a BSDE
satisfies (A1) and (A2).

Let us remark that Lemma 2.1 remains true when assumption (A1) is
replaced by weaker assumptions, for instance those studied in
Bahlali~\cite{B}, Bahlali, Essaky, Hassani and Pardoux~\cite{BEHP}
or Pardoux and Peng~\cite{PP}. However, here, for the sake of
simplicity of the calculation we prefer to work with the Lipschitz
assumption.

   We also shall recall the following both basic results on BSDEs.
   We begin with the well-known comparison theorem (see Theorem 2.2
   in El Karoui, Peng and Quenez~\cite{ElPeQu}).

\bl (Comparison Theorem) Given two coefficients $g_1$ and $g_2$
satisfying (A1) and (A2) and two terminal values $ \xi_1,\ \xi_2 \in
L^{2}(\Omega, {\cal{F}}_{T}, P)$, we denote by $(y^1,z^1)$\ and
$(y^2,z^2)$\ the solution of the BSDE with data $(\xi_1,g_1 )$\ and
$(\xi_2,g_2 )$, respectively. Then we have:

{\rm (i) }(Monotonicity) If  $ \xi_1 \geq \xi_2$  and $ g_1 \geq
g_2, \ a.s.$, then $y^1_t\geq y^2_t,\ a.s.$, for all $t\in [0,
T].$

{\rm (ii)}(Strict Monotonicity) If, in addition to {\rm (i)}, we
also assume that $P(\xi_1 > \xi_2)> 0$, then $P\{y^1_t> y^2_t\}>0,
\ 0 \leq t \leq T,$\ and in particular, $ y^1_0> y^2_0.$ \el

Using the notation introduced in Lemma 2.2 we now suppose that, for
some $g: \Omega\times[0, T]\times{\mathbb{R}}
\times{\mathbb{R}}^{d}\longrightarrow {\mathbb{R}}$\ satisfying (A1)
and (A2), the drivers $g_i, \ i=1, 2,$\ are of the form
$$g_i(s, y_s^i, z_s^i)=g(s, y_s^i, z_s^i)+\varphi_i(s),\ \ \mbox{dsdP-a.e.},$$
where $\varphi_i\in {\cal{H}}_{\mathbb{F}}^{2}(0,T;{\mathbb{R}}).$\
Then, for terminal values $\xi_1,\ \xi_2\ \mbox{belonging to}\
L^{2}(\Omega, {\cal{F}}_{T}, P)$\ we have the following

 \bl The difference of the solutions $(y^1, z^1)$ and $(y^2, z^2)$ of the BSDE with data
 $(\xi_1, g_1)$\ and $(\xi_2, g_2)$, respectively, satisfies
 the following estimate:
 $$
  \begin{array}{ll}
  &|y^1_t-y^2_t|^2+\frac{1}{2}E[\int^T_te^{\beta(s-t)}[|
  y^1_s-y^2_s|^2+ |
  z^1_s-z^2_s|^2]ds|{\cal{F}}_t]
     \\
  \leq& E[e^{\beta(T-t)}|\xi_1-\xi_2|^2|{\cal{F}}_t]+ E[\int^T_te^{\beta(s-t)}
           |\varphi_1(s)-\varphi_2(s)|^2ds|{\cal{F}}_t],\ \mbox{P-a.s.,
           \ for all}\ 0\leq t\leq T,
  \end{array}
  $$
where $\beta=16(1+C^2)$.
  \el
For the proof the reader is referred to Proposition 2.1 in El
Karoui, Peng and Quenez~\cite{ElPeQu} or Theorem 2.3 in
Peng~\cite{Pe1}.

\section{\large{Mean-Field Backward Stochastic Differential
Equations}}

This section is devoted to the study of a new type of BSDEs, the so
called Mean-Field BSDEs.

Let $(\bar{\Omega}, \bar{{\cal{F}}}, \bar{P})=(\Omega\times \Omega,
{\cal{F}}\otimes {\cal{F}}, P\otimes P )$\ be the (non-completed)
product of $(\Omega, {\cal{F}}, P)$\ with itself. We endow this
product space with the filtration
$\bar{{\mathbb{F}}}=\{\bar{{\mathcal{F}}}_t={\cal{F}}\otimes
{\cal{F}}_t,\ 0\leq t \leq T\}$. A random variable $\xi\in
L^0(\Omega, {\cal{F}}, P; {\mathbb{R}}^n)$\ originally defined on
$\Omega$\ is extended canonically to $\bar{\Omega}$:
$\xi'(\omega',\omega)=\xi(\omega'),\ \ (\omega', \omega)\in
\bar{\Omega}=\Omega\times \Omega$. For any $\theta\in
L^1(\bar{\Omega}, \bar{{\cal{F}}}, \bar{P})$\ the variable
$\theta(., \omega): \Omega\rightarrow {\mathbb{R}}$\ belongs to
$L^1({\Omega}, {\cal{F}},P),\ \ P(d\omega)\mbox{-a.s.};$\ we denote
its expectation by
$$E'[\theta(., \omega)]=\int_{\Omega}\theta(\omega', \omega)P(d\omega').$$
Notice that $E'[\theta]=E'[\theta(., \omega)]\in L^1({\Omega},
{\cal{F}},P),$\ and
$$\bar{E}[\theta](=\int_{\bar{\Omega}}\theta d\bar{P}=\int_\Omega E'[\theta(.,
\omega)]P(d\omega))=E[E'[\theta]].$$

\noindent The driver of our mean-field BSDE is a function
$f=f(\omega',\omega, t, y',z', y,z): \bar{\Omega}\times[0,T]\times
{\mathbb{R}} \times {\mathbb{R}}^{d}\times {\mathbb{R}} \times
{\mathbb{R}}^{d} \rightarrow {\mathbb{R}}$\ which is
$\bar{{\mathbb{F}}}$-progressively measurable, for all
$(y',z',y,z)$, and satisfies the following assumptions:
 \vskip0.2cm

(A3) There exists a constant $C\ge 0$  such that, $\bar{P}$-a.s.,
for all $t\in [0, T],\ y_{1}, y_{2}, y'_{1}, y'_{2}\in
{\mathbb{R}},\ z_{1}, z_{2}, z'_{1}, z'_{2}\in {\mathbb{R}}^d,\\
\mbox{ }\hskip1.5cm |f(t, y'_{1}, z'_{1}, y_{1}, z_{1}) - f(t,
y'_{2}, z'_{2}, y_{2}, z_{2})|\leq C(|y'_{1}-y'_{2}|+
|z'_{1}-z'_{2}|+|y_{1}-y_{2}|+ |z_{1}-z_{2}|).$
 \vskip0.2cm

(A4) $f(\cdot,0,0,0,0)\in
{\cal{H}}_{\bar{{\mathbb{F}}}}^{2}(0,T;{\mathbb{R}})$. \vskip0.2cm

\br Let $\beta: \Omega\times [0, T]\rightarrow {\mathbb{R}},\
\gamma: \Omega\times [0, T]\rightarrow {\mathbb{R}}^d$\ be two
square integrable, jointly measurable processes. Then, for our
driver, we can define, for all $(y, z)\in {\mathbb{R}}\times
{\mathbb{R}}^d,\ \mbox{dt}P(d\omega)\mbox{-a.e.},$
$$\begin{array}{rcl} f^{\beta, \gamma}(\omega, t,y,z)&=&E'[f(., \omega, t, \beta'_t, \gamma'_t, y,z)]\\
&=& \int_\Omega f(\omega', \omega, t, \beta_t(\omega'),
\gamma_t(\omega'), y,z)P(d\omega').\end{array}
$$
Indeed, we remark that, for all $(y, z)$,\ due to our assumptions on
the driver f, $(f(., t, \beta'_t, \gamma'_t, y,z))\in
{\cal{H}}_{\bar{{\mathbb{F}}}}^{2}(0,T;{\mathbb{R}}),$\ and thus
$f^{\beta, \gamma}(.,.,y,z)\in {\cal{H}}_{
{\mathbb{F}}}^{2}(0,T;{\mathbb{R}})$. Moreover, with the constant
$C$\ of assumption (A3), for all $(y_1, z_1),\  (y_2, z_2)\in
{\mathbb{R}}\times {\mathbb{R}}^d,\
\mbox{dtP}(d\omega)\mbox{-a.e.},$ $$|f^{\beta, \gamma}(\omega,t,y_1,
z_1)-f^{\beta, \gamma}(\omega,t,y_2, z_2)|\leq
C(|y_1-y_2|+|z_1-z_2|).$$ Consequently, there is an
${\mathbb{F}}$-progressively measurable version of $f^{\beta,
\gamma}(.,.,y, z),\ (y,z)\in {\mathbb{R}}\times {\mathbb{R}}^d$,
such that $f^{\beta, \gamma}(\omega,t,., .)$\ is
$\mbox{dtP}(d\omega)\mbox{-a.e.}$\ defined and Lipschitz in $(y,z)$;
its Lipschitz constant is that introduced in (A3). \er

We now can state the main result of this section.

\bt Under the assumptions (A3) and (A4), for any random variable
$\xi\in L^2(\O, {\cal{F}}_T,$ $P),$ the mean-field BSDE
 \be Y_t = \xi + \int_t^TE'[f(s,Y'_s,Z'_s,Y_s,Z_s)]ds - \int^T_tZ_s\,
dB_s,\q 0\le t\le T, \label{BSDE1} \ee
 has a unique adapted solution
$$(Y_t, Z_t)_{t\in [0, T]}\in {\cal{S}}_{{\mathbb{F}}}^2(0, T; {\mathbb{R}})\times
{\cal{H}}_{{\mathbb{F}}}^{2}(0,T;{\mathbb{R}}^{d}). $$ \et \br We
emphasize that, due to our notations, the driving coefficient of
(3.1) has to be interpreted as follows
$$\begin{array}{rcl} E'[f(s,Y'_s, Z'_s, Y_s, Z_s)](\omega)&=& E'[f(s,Y'_s, Z'_s,
Y_s(\omega), Z_s(\omega)]\\
&=& \int_\Omega f(\omega', \omega, s, Y_s(\omega'),
Z_s(\omega'),Y_s(\omega), Z_s(\omega))P(d\omega').\end{array}
$$
\er

\noindent {\bf Proof}. We first introduce a norm on the space
${\cal{H}}_{{\mathbb{F}}}^{2}(0,T;{\mathbb{R}}\times{\mathbb{R}}^{d})$\
which is equivalent to the canonical norm:
$$||v(\cdot)||_\beta=\left\{E\int_0^T|v_s|^2e^{\beta s}ds\right\}^{\frac{1}{2}}, \ \ \ \  \ \beta>0.$$
The parameter $\beta$\ will be specified later.

\noindent \underline{Step 1}: For any $(y,z)\in
{\cal{H}}_{{\mathbb{F}}}^{2}(0,T;{\mathbb{R}}\times{\mathbb{R}}^{d})$\
there exists a unique solution $(Y, Z)\in
{\cal{S}}_{{\mathbb{F}}}^2(0, T; {\mathbb{R}})\times
{\cal{H}}_{{\mathbb{F}}}^{2}(0,T;{\mathbb{R}}^{d})$\ to the
following BSDE: \be Y_t = \xi + \int_t^TE'[f(s,y'_s,z'_s,Y_s,Z_s)]ds
- \int^T_tZ_s\, dB_s,\q 0\le t\le T. \label{BSDE2} \ee

Indeed, we define $g^{(y,z)}(s,\mu,\nu)=
E'[f(s,y'_s,z'_s,\mu,\nu)]$. Then, due to Remark 3.1,
 $g^{(y,z)}(s,\mu,\nu)$\ satisfies (A1) and (A2), and from Lemma 2.1 we
know there exists a unique solution $(Y, Z)\in
{\cal{S}}_{{\mathbb{F}}}^2(0, T; {\mathbb{R}})\times
{\cal{H}}_{{\mathbb{F}}}^{2}(0,T;{\mathbb{R}}^{d})$\ to the BSDE
(3.2).

\vskip 0.3cm \noindent \underline{Step 2}: The result of Step 1
allows to introduce the mapping $(Y., Z.)=I[(y., z.)]:
{\cal{H}}_{{\mathbb{F}}}^{2}(0,T;{\mathbb{R}}\times{\mathbb{R}}^{d})\rightarrow
{\cal{H}}_{{\mathbb{F}}}^{2}(0,T;{\mathbb{R}}\times{\mathbb{R}}^{d})$\
by the equation\be Y_t = \xi + \int_t^TE'[f(s,y'_s,z'_s,Y_s,Z_s)]ds
- \int^T_tZ_s\, dB_s,\q 0\le t\le T. \label{BSDE3} \ee For any
$(y^1,z^1),\ (y^2,z^2)\in
{\cal{H}}_{{\mathbb{F}}}^{2}(0,T;{\mathbb{R}}\times{\mathbb{R}}^{d})$\
we put $(Y^1, Z^1)=I[(y^1,z^1)],\ (Y^2, Z^2)=I[(y^2,z^2)]$,
$(\hat{y},\hat{z})=(y^1-y^2, z^1-z^2)$\ and
$(\hat{Y},\hat{Z})=(Y^1-Y^2, Z^1-Z^2).$\ Then, by applying It\^{o}'s
formula to $e^{\beta s}|\hat{Y}_s|^2$\ and by using that $Y^1,\ \
Y^2\in {\cal{S}}_{{\mathbb{F}}}^2(0, T; {\mathbb{R}})$\ we get
$$\begin{array}{rcl}
 & & |\hat{Y}_t|^2 +E[\int_t^T e^{\beta(r-t)}\beta|\hat{Y}_r|^2dr|{\cal{F}}_t]+E[\int_t^T
 e^{\beta(r-t)}|\hat{Z}_r|^2dr|{\cal{F}}_t]\\
 &= & E[\int_t^T e^{\beta(r-t)}2\hat{Y}_r(g^{(y^1,z^1)}(r, Y^1_r, Z^1_r)-g^{(y^2,z^2)}(r, Y^2_r,
 Z^2_r))dr|{\cal{F}}_t],\ \ t\in [0, T].\\
 \end{array}$$

\noindent From assumption (A3) we obtain
$$\begin{array}{rcl}
 & & (\frac{\beta}{2}-2C-2C^2)E[\int_0^T e^{\beta r} |\hat{Y}_r|^2dr]+\frac{1}{2}E[\int_0^T
 e^{\beta r}|\hat{Z}_r|^2dr]\\
 &\leq& \frac{4C^2}{\beta}\{E[\int_0^T e^{\beta r}|\hat{y}_r|^2dr]+E[\int_0^T e^{\beta r}|\hat{z}_r|^2dr]\}.\\
 \end{array}$$
Thus, taking $\beta=16C^2+4C+1$\ we get
$$E[\int_0^Te^{\beta r}(|\hat{Y}_r|^2+|\hat{Z}_r|^2)dr]\leq \frac{1}{2}
E[\int_0^Te^{\beta r}(|\hat{y}_r|^2+|\hat{z}_r|^2)dr],$$ that is,
$||(\hat{Y}, \hat{Z})||_\beta \leq \frac{1}{\sqrt{2}}||(\hat{y},
\hat{z})||_\beta.$\ Consequently, $I$\ is a contraction on
${\cal{H}}_{{\mathbb{F}}}^{2}(0,T;{\mathbb{R}}\times{\mathbb{R}}^{d})$\
endowed with the norm $||.||_\beta$, and from the contraction
mapping theorem we know that there is a unique fixed point $(Y,\
Z)\in
{\cal{H}}_{{\mathbb{F}}}^{2}(0,T;{\mathbb{R}}\times{\mathbb{R}}^{d})$\
such that $I(Y, Z)=(Y, Z).$\ On the other hand, from Step 1 we
already know that if $I(Y, Z)=(Y, Z)$\ then $(Y,\ Z)\in
{\cal{S}}_{{\mathbb{F}}}^2(0, T; {\mathbb{R}})\times
{\cal{H}}_{{\mathbb{F}}}^{2}(0,T;{\mathbb{R}}^{d}).$\endpf

Using the notation introduced in Theorem 3.1 we now suppose that,
for some $f: \bar{\Omega}\times[0, T]\times{\mathbb{R}}
\times{\mathbb{R}}^{d}\times{\mathbb{R}}
\times{\mathbb{R}}^{d}\rightarrow {\mathbb{R}}$\ satisfying (A3) and
(A4), the drivers $f_i, \ i=1, 2,$\ are of the form
$$f_i(s, Y_s^{i'}, Z_s^{i'}, Y_s^i, Z_s^i)=f(s, Y_s^{i'}, Z_s^{i'},
Y_s^i, Z_s^i)+\varphi_i(s),\ \ \mbox{dsd}\bar{P}\mbox{-a.e.},\ i=1, 2,$$
where $\varphi_i\in
{\cal{H}}_{\bar{{\mathbb{F}}}}^{2}(0,T;{\mathbb{R}})$\ and $(Y^i,
Z^i)$\ is the solution of Mean-Field BSDE (3.1) with data $(f_i,
\xi_i),\ \ i=1, 2.$\ Then, for arbitrary terminal values $\xi_1,\
\xi_2\ \mbox{belonging to}\ L^{2}(\Omega, {\cal{F}}_{T}, P)$\ we
have the following

 \bl The difference of the solutions $(Y^1, Z^1)$ and $(Y^2, Z^2)$ of
 BSDE (3.1) with the data $(\xi_1, f_1)$\ and $(\xi_2, f_2)$, respectively,
 satisfies the following estimate:
 $$
  \begin{array}{ll}
  &E[|Y^1_t-Y^2_t|^2]+\frac{1}{2}E[\int^T_te^{\beta(s-t)}(|
  Y^1_s-Y^2_s|^2+ |
  Z^1_s-Z^2_s|^2)ds]
     \\
  \leq& E[e^{\beta(T-t)}|\xi_1-\xi_2|^2]+ \bar{E}[\int^T_te^{\beta(s-t)}
           |\varphi_1(s)-\varphi_2(s)|^2ds],\ \mbox{for all}\ 0\leq t\leq T,
  \end{array}
  $$
where $\beta=16(1+C^2)$.
  \el
The proof uses a similar argument as that of the proof of Theorem
3.1 and is therefore omitted.

Now we discuss the comparison principle for Mean-Field BSDE. We
first give two counterexamples to show that if the driver $f$\
depends on $z'$\ or $f$\ is decreasing with respect to $y'$\ we
can't get the comparison theorem.

\noindent {\bf Example 3.1}\ \ For $d=1$\ we consider the Mean-Field
BSDE (3.1) with time horizon $T=1$,\ with driver
$f(\omega',\omega,s, y', z', y, z)=-z'$\ and two different terminal
values $\xi_1,\ \xi_2\in L^{2}(\Omega, {\cal{F}}_{T}, P).$\ Let us
denote the associated solutions by $(Y^1, Z^1)$\ and $(Y^2, Z^2)$,
respectively. Then, \be Y^i_t = \xi_i + \int_t^TE[-Z^i_s]ds -
\int^T_tZ^i_sdB_s,\q \ 0\le t\le 1,\ i=1,\ 2.\ee  We let
$\xi_1=-(B_1^+)^3$\ and define $\widetilde{Y}^1_t= Y^1_t+
\int_t^TE[Z^1_s]ds.$\ Then $(\widetilde{Y}^1, Z^1)$\ is the unique
solution of the BSDE $\widetilde{Y}^1_t=\xi_1 - \int^T_tZ^1_sdB_s,\
t\in [0, T].$\ Thus,
 $\widetilde{Y}^1_0=E[\xi_1]=-E[(B_1^+)^3]=-\frac{2}{\sqrt{2\pi}}.$\ On the
other hand, from the Clark-Ocone formula we know that $Z^1$\ is the
predictable projection of the Malliavin derivative $(D_t\xi)_{t\in
[0, T]}$\ of $\xi$\ $(D_t\xi\ \mbox{denotes the Malliavin derivative
of}\ \ \xi\ \mbox{at}\ \mbox{time}\ t;$\ the interested reader is
referred, e.g., to Nualart~\cite{N}). This implies
$E[Z^1_t]=E[D_t\xi]=E[-3(B_1^+)^2]=-\frac{3}{2}, \ t\in [0, 1].$
 Therefore,
$Y^1_0=\widetilde{Y}^1_0-\int_0^TE[Z^1_s]ds=-\frac{2}{\sqrt{2\pi}}+\frac{3}{2}>0.$\
Let now $\xi_2=0$. Then, obviously, $(Y^2, Z^2)=(0, 0).$\ Hence, we
have $Y^1_0>Y^2_0$\ although $\xi_1\leq \xi_2,\ \ \mbox{P-a.s.}$\
and $P\{\xi_1<\xi_2\}>0.$

\noindent {\bf Example 3.2}\ \ Let again $d=1$. We consider
Mean-Field BSDE (3.1) driven by the function $f(\omega',\omega,s,
y', z', y, z)=-y',$\ with time horizon $T=2$\ and two different
terminal values $\xi_1,\ \xi_2\in L^{2}(\Omega, {\cal{F}}_{T}, P):$
\be Y^i_t = \xi_i + \int_t^TE[-Y^i_s]ds - \int^T_tZ^i_sdB_s,\q  0\le
t\le 2,\ i=1,\ 2.\ee  By choosing first $\xi_1=(B_1)^2$\ we have
$E[Y^1_t]=e^{-(2-t)},\ \ t\in [0, 2].$\ Furthermore, for $t\in [1,
2]$\ we have
$Y^1_t=(B_1)^2-\int_t^2e^{-(2-s)}ds=(B_1)^2-(1-e^{-(2-t)})\
\mbox{and}\ Z^1_t=0,\ \mbox{P-a.s.} $\ Consequently,
$Y^1_1=(B_1)^2-(1-e^{-1})<0$\ on the set $\{(B_1)^2 < 1-e^{-1}\}$\
which is of strictly positive probability. Finally, for $\xi_2=0$\
we have the solution $(Y^2, Z^2)=(0, 0).$\ Therefore, in our example
$P(Y^1_1<Y^2_1)>0,$\ although $\xi_1> \xi_2,\ \ \mbox{P-a.s.}$

The above both examples show that we cannot hope to have the
comparison principle between two Mean-Field BSDEs whose drivers
depend both on $z'$\ or are both decreasing in $y'$.

\bt (Comparison Theorem) Let $f_i=f_i(\bar{\omega}, t, y', z', y,
z),\ i=1, 2,$\ be two drivers satisfying the standard assumptions
(A3) and (A4). Moreover, we suppose

\noindent{\rm(i)} One of both coefficients is independent of $z'$;

\noindent{\rm(ii)} One of both coefficients is nondecreasing in
$y'$.

Let $\xi_1,\ \xi_2\in L^{2}(\Omega, {\cal{F}}_{T}, P)$\ and denote
by $(Y^1, Z^1)$\ and $(Y^2, Z^2)$\ the solution of the mean-field
BSDE (3.1) with data $(\xi_1, f_1)$\ and $(\xi_2, f_2)$,
respectively. Then if  $\xi_1\leq \xi_2,\ \mbox{P-a.s.}$, and
$f_1\leq f_2,\ \bar{P}\mbox{-a.s.},$\ it holds that also $Y_t^1\leq
Y_t^2,\ t\in [0, T],\ \mbox{P-a.s.}$\et

\br The conditions {\rm(i)} and {\rm(ii)} of Theorem 3.2 are, in
particular, satisfied, if they hold for the same driver $f_j$\ but
also if {\rm(i)} is satisfied by one driver and {\rm(ii)} by the
other one. \er

\noindent {\bf Proof.} Without loss of generality, we assume that
{\rm(i)} is satisfied by $f_1$ and {\rm(ii)} by $f_2$. Then, with
the notation $(\bar{Y}, \bar{Z}):= (Y^1-Y^2, Z^1-Z^2),\
\bar{\xi}:=\xi_1-\xi_2,$ \be \bar{Y}_t = \bar{\xi} +
\int_t^TE'[f_1(s,Y^{1'}_s, Y^1_s,Z^1_s)-f_2(s,Y^{2'}_s, Z^{2'}_s,
Y^2_s,Z^2_s)]ds - \int^T_t\bar{Z}_s\, dB_s,\q 0\le t\le T,
\label{BSDE1} \ee and from It\^{o}'s formula applied to an
appropriate approximation of $\varphi(y)=(y^+)^2,\ y\in
{\mathbb{R}},$\ in which we take after the limit, we obtain
$$\begin{array}{rcl}
& & E[(\bar{Y}_t^+)^2]+E[\int_t^TI_{\{\bar{Y}_s>0\}}|\bar{Z}_s|^2ds]\\
&=& 2E[\int_t^T\bar{Y}_s^+(E'[f_1(s,Y^{1'}_s,
Y^1_s,Z^1_s)-f_2(s,Y^{2'}_s, Z^{2'}_s, Y^2_s,Z^2_s)])ds]\\
&\leq & 2E[\int_t^T\bar{Y}_s^+(E'[f_1(s,Y^{1'}_s,
Y^1_s,Z^1_s)-f_1(s,Y^{1'}_s, Y^2_s,Z^2_s)])ds]\\
& & + 2E[\int_t^T\bar{Y}_s^+(E'[f_2(s,Y^{1'}_s,Z^{2'}_s,
Y^2_s,Z^2_s)-f_2(s,Y^{2'}_s, Z^{2'}_s, Y^2_s,Z^2_s)])ds],\ \ t\in[0, T].\\
\end{array}
$$
Since $f_2(s, y',z', y, z)$\ is nondecreasing in $y'$\ we get from
(A3)
$$\begin{array}{rcl}
& & E[\bar{Y}_s^+(E'[f_1(s,Y^{1'}_s,
Y^1_s,Z^1_s)-f_1(s,Y^{1'}_s, Y^2_s,Z^2_s)])]\\
& & + E[\bar{Y}_s^+(E'[f_2(s,Y^{1'}_s,Z^{2'}_s,
Y^2_s,Z^2_s)-f_2(s,Y^{2'}_s,Z^{2'}_s, Y^2_s,Z^2_s)])]\\
&\leq &
E[\bar{Y}_s^+(E'[C(Y^{1'}_s-Y^{2'}_s)^++C|\bar{Y}_s|+C|\bar{Z}_s|])]\\
&= &
C(E[\bar{Y}_s^+])^2+CE[(\bar{Y}_s^+)^2]+CE[\bar{Y}_s^+|\bar{Z}_s|]\\
&\leq &
(2C+C^2)E[(\bar{Y}_s^+)^2]+\frac{1}{4}E[I_{\{\bar{Y}_s>0\}}|\bar{Z}_s|^2], \ \ s\in [0, T].\\
\end{array}
$$
Consequently,
$$E[(\bar{Y}_t^+)^2]+\frac{1}{2}E[\int_t^TI_{\{\bar{Y}_s>0\}}|\bar{Z}_s|^2ds]\leq
(4C+2C^2)\int_t^TE[(\bar{Y}_s^+)^2]ds,\ \ t\in [0, T],
$$
from where we can conclude with the help of Gronwall's Lemma that
$Y^1_t- Y^2_t=\bar{Y}_t\leq 0,\ t\in [0, T],\ \mbox{P-a.s.}$\endpf

We also have a converse comparison theorem.

\bt (Converse Comparison Theorem) We retake the assumptions of
Theorem 3.2 and suppose that, additionally, for some $t\in [0, T],\
Y^1_t=Y^2_t,\ \mbox{P-a.s.}$\ Then

\noindent{\rm(i)} $Y^1_s=Y^2_s,\ s\in [t, T],\ \mbox{P-a.s.}$, and

\noindent{\rm(ii)} If $f_2$\ satisfies {\rm(ii)} of Theorem 3.2 then
$E'[f_1(s,Y^{1'}_s, Z^{1'}_s, Y^2_s,Z^2_s)]=E'[f_2(s,Y^{2'}_s,
Z^{2'}_s, Y^2_s,Z^2_s)],$\ $\mbox{dsdP-a.e.}$\ on $[t, T]$, and if
$f_1$\ satisfies {\rm(ii)} then the symmetric result holds.\et

\noindent {\bf Proof.} We use the notation introduced in the proof
of Theorem 3.2 and suppose again that $f_1$\ satisfies {\rm(i)} and
$f_2$\ {\rm(ii)}. Then, from the Lipschitz property of $f_1$, there
exist some $\bar{\mathbb{F}}$-progressively measurable, bounded
processes $\alpha,\ \beta$, defined over $\bar{\Omega}\times [0,
T]$, such that
$$f_1(s, Y_s^{1'}, Y_s^{1}, Z_s^{1})-f_1(s, Y_s^{1'}, Y_s^{2}, Z_s^{2})
=\alpha_s\bar{Y}_s+\beta_s\bar{Z}_s,\ s\in [0, T],$$
and since $Y^{1},\ Z^{1},\ Y^{2},\ Z^{2}$\ don't depend on
$\omega'$,
$$E'[f_1(s, Y_s^{1'}, Y_s^{1}, Z_s^{1})-f_1(s, Y_s^{1'}, Y_s^{2}, Z_s^{2})]
=E'[\alpha_s]\bar{Y}_s+E'[\beta_s]\bar{Z}_s,\ s\in [0, T].$$
Thus, from It\^{o}'s formula,
$$\begin{array}{rcl}
\bar{Y}_t&=&\exp\{\int_t^TE'[\alpha_r]dr\}\bar{\xi}+\int_t^T\exp\{\int_t^s
E'[\alpha_r]dr\}
E'[f_1(s, Y_s^{1'}, Y_s^{2}, Z_s^{2})-f_2(s, Y_s^{2'},
Z_s^{2'},Y_s^{2},
Z_s^{2})]ds\\
& &-\int_t^T\exp\{\int_t^sE'[\alpha_r]dr\}\bar{Z}_sd\tilde{B}_s,\ \
\mbox{P-a.s.,}\end{array}$$ where
$\tilde{B}_s=B_s-\int_0^sE'[\beta_r]dr,\ s\in [0, T].$\ It is well
known that $\tilde{B}=(\tilde{B}_s)$\ is an $({\mathbb{F}},
\tilde{P})$-Brownian motion with
$\tilde{P}:=\exp\{\int_0^TE'[\beta_s]dB_s-\frac{1}{2}\int_0^T|E'[\beta_s]|^2ds\}P.$\
From the boundedness of $\alpha$\ and $\beta$\ we then deduce easily
that $\int_t^T\exp\{\int_t^sE'[\alpha_r]dr\}\bar{Z}_sd\tilde{B}_s$\
is an $({\mathbb{F}}, \tilde{P})$-martingale increment.
Consequently,
$$\begin{array}{rcl}
0=\bar{Y}_t& = &E_{\tilde{P}}[\exp\{\int_t^uE'[\alpha_r]dr\}\bar{Y}_u|{\cal{F}}_t]\\
& &+E_{\tilde{P}}[\int_t^u\exp\{\int_t^sE'[\alpha_r]dr\} E'[f_1(s,
Y_s^{1'}, Y_s^{2}, Z_s^{2})-f_2(s, Y_s^{2'}, Z_s^{2'},Y_s^{2},
Z_s^{2})]ds|{\cal{F}}_t],\\
\ & &\hskip9cm \ \mbox{P-a.s., for all}\ t\leq u\leq T.\end{array}$$
To conclude it suffices now to recall that, due to Theorem 3.2,
$\bar{Y}_u=Y_u^1-Y_u^2\leq 0,\ \mbox{P-a.s.}$\ and
$$f_1(s,Y_s^{1'}, Y_s^{2}, Z_s^{2})\leq f_2(s,Y_s^{1'},Z_s^{2'}, Y_s^{2},
Z_s^{2})\leq f_2(s,Y_s^{2'}, Z_s^{2'}, Y_s^{2}, Z_s^{2}),\
\mbox{dsd}\bar{P}\mbox{-a.e. on}\ \bar{\Omega}\times[t, T].$$
\endpf
\section{\large{Mean-Field Stochastic Differential
Equations }}

\hskip1cm We consider measurable functions $b:\bar{\Omega}\times
[0,T]\times {\mathbb{R}}^n\times {\mathbb{R}}^n\rightarrow
{\mathbb{R}}^n \ $ and
         $\sigma:\bar{\Omega}\times[0,T]\times  {\mathbb{R}}^n\times
         {\mathbb{R}}^n\rightarrow {\mathbb{R}}^{n\times d}$
which are supposed to satisfy the following conditions:
 $$
  \begin{array}{ll}
\mbox{(i)}&b(\cdot,0,0)\ \mbox{and}\ \sigma(\cdot,0, 0)\ \mbox{are}
\ \bar{{\mathbb{F}}}-\mbox{progressively measurable continuous
processes and
there exists}\\ & \mbox{some constant}\ C>0\  \mbox{such that}\\
 &\hskip 1cm|b(t,x',x)|+|\sigma(t,x',x)|\leq C(1+|x|), a.s.,\
                                  \mbox{for all}\ 0\leq t\leq T,\ x, x'\in {\mathbb{R}}^n;\\
\mbox{(ii)}&b\ \mbox{and}\ \sigma\ \mbox{are Lipschitz in}\ x,\ x',\ \mbox{i.e., there is some constant}\ C>0\ \mbox{such that}\\
           &\hskip 1cm|b(t,x'_1,x_1)-b(t, x'_2,x_2)|+|\sigma(t,x'_1,x_1)-\sigma(t, x'_2 ,x_2)|\leq C(| x'_1-x'_2|+| x_1-x_2|),\ a.s.,\\
 & \hbox{ \ \ }\hskip7cm\mbox{for all}\ 0\leq t \leq T,\ x_1,x'_1, \ x_2, x'_2\in {\mathbb{R}}^n.\\
 \end{array}
  \eqno{\mbox{(H4.1)}}
  $$\par
  We now study the following SDE parameterized by the
  initial condition $(t,\zeta)\in[0,T]\times L^2(\Omega,{\cal{F}}_t,P;{\mathbb{R}}^n)$:
  \be
  \left\{
  \begin{array}{rcl}
  dX_s^{t,\zeta}&=&E'[b(s,(X_s^{t,\zeta})',X_s^{t,\zeta})]ds+E'[\sigma(s,(X_s^{t,\zeta})',X_s^{t,\zeta})]dB_s,\ s\in[t,T],\\
  X_t^{t,\zeta}&=&\zeta.
  \end{array}
  \right.
  \ee
Recall that
$$E'[b(s,(X_s^{t,\zeta})',X_s^{t,\zeta})](\omega)=
\int_\Omega
b(\omega',\omega,s,X_s^{t,\zeta}(\omega'),X_s^{t,\zeta}(\omega))P(d\omega'),\
\ \ \omega\in \Omega.$$ \bt Under the assumption (H4.1) SDE (4.1)
has a unique strong solution. \et

\noindent The proof is made in two steps like that of Theorem 3.1
and it uses standard arguments for forward SDEs. Since the proof is
straight-forward we prefer to omit it.

\br From standard arguments we also can get that, for any $p\geq
2,$\ there exists $C_{p}\in {\mathbb{R}}$\ such that, for all
$t\in[0,T]\ \mbox{and}\ \zeta,\zeta'\in
L^p(\Omega,{\cal{F}}_t,P;{\mathbb{R}}^n),$
 \be
 \begin{array}{rcl}
 E[\sup\limits_{t\leq s\leq T}| X_s^{t,\zeta}-X_s^{t,\zeta'}|^p|{\cal{F}}_t]
                             &\leq& C_{p}|\zeta-\zeta'|^p, \ \ a.s.,\\
  E[\sup\limits_{t\leq s\leq T}| X_s^{t,\zeta}|^p|{\cal{F}}_t]
                       &\leq& C_{p}(1+|\zeta|^p),\ \  a.s.,\\
  E[\sup\limits_{t\leq s\leq t+\delta}|X^{t,\zeta}_{s} -\zeta|^p|{{\mathcal{F}}_t}]
  &\leq &
C_{p}(1+|\zeta|^p)\delta^{\frac{p}{2}},
 \end{array}
\ee
 $\mbox{P-a.s.},\ \mbox{for all}\ \delta>0$\ with $t+\delta\leq
T.$

These in the classical case well-known standard estimates can be
consulted, for instance, in Ikeda, Watanabe~\cite{IW}, pp.166-168
and also in Karatzas, Shreve~\cite{KSH}, pp.289-290. We also
emphasize that the constant $C_{p}$ in (4.2) only depends on the
Lipschitz and the growth constants of $b$ and $\sigma$. \er

\section{\large{Decoupled Mean-Field Forward-Backward SDE and Related DPP}}

\hskip1cm In this section we study a decoupled Mean-Field
forward-backward SDE and its relation with PDEs. Given continuous
functions $b:\bar{\Omega}\times[0,T]\times  {\mathbb{R}}^n\times
{\mathbb{R}}^n\rightarrow {\mathbb{R}}^n \ $ and
$\sigma:\bar{\Omega}\times[0,T]\times
 {\mathbb{R}}^n\times {\mathbb{R}}^n\rightarrow
{\mathbb{R}}^{n\times d}$ which are supposed to satisfy the
conditions {\rm(i)} and {\rm(ii)} of (H4.1) and an arbitrary $x_0\in
{\mathbb{R}}^n,$\ we consider the following SDE parameterized by the
  initial condition $(t,\zeta)\in[0,T]\times L^2(\Omega,{\cal{F}}_t,P;{\mathbb{R}}^n)$:
  \be
  \left\{
  \begin{array}{rcl}
  dX_s^{t,\zeta}&=&E'[b(s,(X_s^{0,x_0})',X_s^{t,\zeta})]ds+E'[\sigma(s,(X_s^{0,x_0})',X_s^{t,\zeta})]dB_s,\ s\in[t,T],\\
  X_t^{t,\zeta}&=&\zeta.
  \end{array}
  \right.
  \ee
Under the assumption (H4.1), SDE (5.1) has a unique strong solution.
Indeed, from Theorem 4.1 we first deduce the existence and
uniqueness of the process $X^{0,x_0}\in
{\cal{S}}_{{\mathbb{F}}}^2(0, T; {\mathbb{R}}^n)$\ to the Mean-Field
SDE (5.1). Once knowing $X^{0,x_0}$, SDE (5.1) becomes a classical
equation with the coefficients
$\tilde{b}(\omega,s,x)=E'[b(\omega',\omega,s,X_s^{0,x_0}(\omega'),x)]$\
and
$\tilde{\sigma}(\omega,s,x)=E'[\sigma(\omega',\omega,s,X_s^{0,x_0}(\omega'),x)].$\
Combining estimate (4.2) for $(t, \zeta)=(0, x_0)$\ with standard
arguments for SDEs we obtain (4.2) also for equation (5.1).

Let now be given two real-valued functions $f(t, x',x, y',y,z)$ and
$\Phi(x',x)$ which shall satisfy the following conditions:
$$
\begin{array}{ll}
\mbox{(i)}&\Phi:\bar{\Omega}\times {\mathbb{R}}^n\times
{\mathbb{R}}^n\rightarrow {\mathbb{R}} \ \mbox{is an}\
\bar{{\cal{F}}}_T\otimes{\cal{B}}({\mathbb{R}}^n)
             \mbox{-measurable random variable and}\\
          & f:\bar{\Omega}\times [0,T]\times {\mathbb{R}}^n\times {\mathbb{R}}^n\times {\mathbb{R}}\times {\mathbb{R}}\times
          {\mathbb{R}}^d \rightarrow {\mathbb{R}}\ \mbox{is a measurable process}\ \mbox{such that} \\
          &f(\cdot, x',x, y',y,z)\ \mbox{is}\ \bar{{\mathbb{F}}} \mbox{-adapted, for all $(x', x, y',y,z)\in
          {\mathbb{R}}^n\times {\mathbb{R}}^n\times {\mathbb{R}}\times {\mathbb{R}}\times
          {\mathbb{R}}^d $;}\\
\mbox{(ii)}&\mbox{There exists a constant}\ C>0\ \mbox{such that}\\
          &| f(t,x'_1,x_1,y'_1,y_1,z_1)-f(t,x'_2,x_2,y'_2,y_2,z_2)| +| \Phi(x'_1,x_1)-\Phi(x'_2,x_2)|\\
          &\hskip 2cm\leq C(|x'_1-x'_2|+ |x_1-x_2|+|y'_1-y'_2|+|y_1-y_2|+|z_1-z_2|),\ \  a.s.,\\
&\hskip 1cm\mbox{for all}\ 0\leq t\leq T,\ x_1,x'_1, x_2, x'_2\in
{\mathbb{R}}^n,\ y_1,y'_1, y_2, y'_2 \in {\mathbb{R}}\
 \mbox{and}\ z_1, z_2\in {\mathbb{R}}^d;\\
          \end{array}$$
$$
\begin{array}{ll}

 \mbox{(iii)}&f\ \mbox{and}\ \Phi \ \mbox{satisfy a linear growth condition, i.e., there exists some}\ C>0\\\
    & \mbox{such that, a.s.},\ \mbox{for all}\ x',\ x\in
    {\mathbb{R}}^n,\\
    &\hskip 2cm|f(t,x', x,0,0,0)| + |\Phi(x', x)| \leq C(1+|x|+|x'|);\\
\mbox{(iv)}& f(\bar{\omega},t,x',x,y',y,z)\ \mbox{is continuous in}\ t\ \mbox{for all}\ (x',x,y',y,z),\ \mbox{P(d}\bar{\omega})\mbox{-a.s.};\\
\mbox{(v)}& f(t,x', x, y',y, z)\ \mbox{is nondecreasing with respect to}\ y'.\\
\end{array}
\eqno{\mbox{(H5.1)}} $$ We consider the following BSDE:
\be \left\{
\begin{array}{rcl}
-dY_s^{t,\zeta}&=&E'[f(s,(X_s^{0,x_0})',X_s^{t,\zeta},
(Y_s^{0,x_0})',Y_s^{t,\zeta},Z_s^{t,\zeta})]ds
-Z_s^{t,\zeta}dB_s,\ s\in [t, T],\\
Y_T^{t,\zeta}&=&E'[\Phi((X_T^{0,x_0})',X_T^{t,\zeta})].\\
\end{array}
\right. \ee We first consider the equation (5.2) for $(t,\zeta)=(0,
x_0)$: We know from Theorem 3.1 that there exists a unique solution
$(Y^{0, x_0}, Z^{0, x_0})\in {\cal{S}}_{{\mathbb{F}}}^2(0,
T;{\mathbb{R}})\times{\cal{H}}_{{\mathbb{F}}}^{2}(0,T;{\mathbb{R}}^{d})$\
to the Mean-Field BSDE (5.2). Once we have $(Y^{0, x_0}, Z^{0,
x_0}),$\ equation (5.2) becomes a classical BSDE whose coefficients
$\tilde{f}(\omega,
s,X_s^{t,\zeta},y,z)=E'[f(.,\omega,s,(X_s^{0,x_0})',$
$X_s^{t,\zeta}, (Y_s^{0,x_0})',y,z)]$\ satisfies the assumptions
(A1) and (A2), and $\tilde{\Phi}(\omega,
X_T^{t,\zeta}(\omega))=E'[\Phi(.,\omega,(X_T^{0,x_0})',X_T^{t,\zeta})]\in
L^2(\Omega, {\cal{F}}_T, P)$. Thus, from Lemma 2.1 we know that
there exists a unique solution $(Y^{t,\zeta}, Z^{t,\zeta})\in
{\cal{S}}_{{\mathbb{F}}}^2(0, T;
{\mathbb{R}})\times{\cal{H}}_{{\mathbb{F}}}^{2}(0,T;{\mathbb{R}}^{d})$\
to equation (5.2).

By combining classical BSDE estimates (see, e.g., Proposition 4.1 in
Peng~\cite{Pe1}; or Proposition 4.1 in El Karoui, Peng and
Quenez~\cite{ElPeQu}) with the techniques presented above we see
that there exists a constant $C$\ such that, for any $t\in [0, T]$\
and $\zeta,\ \zeta'\in L^2(\Omega,{\cal{F}}_t,P;{\mathbb{R}}^n),$\
we have the following estimates:

$\mbox{}\hskip3cm\mbox{\rm(i)} E[\sup\limits_{t\leq s\leq
T}|Y_s^{t,\zeta}|^2 + \int_t^T|Z_s^{t,\zeta}|^2ds|{{\cal{F}}_t}]\leq
C(1+|\zeta|^2),\
a.s.; $\\

$\mbox{}\hskip3cm\mbox{\rm(ii)}E[\sup\limits_{t\leq s\leq
T}|Y_s^{t,\zeta}-Y_s^{t,\zeta'}|^2+\int_t^T|Z_s^{t,\zeta}-Z_s^{t,\zeta'}|^2ds|{{\cal{F}}_t}]\leq
C|\zeta-\zeta'|^2,\  a.s. $\\ In particular, \be
 \begin{array}{lll}
\mbox{\rm(iii)}&|Y_t^{t,\zeta}|\leq C(1+|\zeta|),\  a.s.; \hskip3cm\\
\mbox{\rm(iv)}&|Y_t^{t,\zeta}-Y_t^{t,\zeta'}|\leq C|\zeta-\zeta'|,\  a.s.\hskip3cm \\
\end{array}
\ee Here the constant $C>0$\ depends only on the Lipschitz and the
growth constants of $b$,\ $\sigma$, $f$\ and $\Phi$.
\vskip 0.3cm
 Let us now introduce the random field:
\be u(t,x)=Y_s^{t,x}|_{s=t},\ (t, x)\in [0, T]\times{\mathbb{R}}^n,
\ee where $Y^{t,x}$ is the solution of BSDE (5.2) with $x \in
{\mathbb{R}}^n$\ at the place of $\zeta\in
L^2(\Omega,{\cal{F}}_t,P;{\mathbb{R}}^n).$\\

As a consequence of (5.3) we immediately have that, for all $t \in
[0, T] $, P-a.s.,
 \be
\begin{array}{ll}
\mbox{(i)}&| u(t,x)-u(t,y)| \leq C|x-y|,\ \mbox{for all}\ x, y\in {\mathbb{R}}^n;\\
\mbox{(ii)}&| u(t,x)|\leq C(1+|x|),\ \mbox{for all}\ x\in {\mathbb{R}}^n.\\
\end{array}
\ee
 \br In the general situation $u$ is an adapted random function, that is, for any $x\in
 {\mathbb{R}}^n,\ u(\cdot,x)$ is an ${\mathbb{F}}-$adapted real-valued process. Indeed, recall that
  $b, \sigma \ \mbox{and}\ f$\ all are $\bar{{\mathbb{F}}}$-adapted random functions while $\Phi$\ is $\bar{{\cal{F}}}_T$-measurable.
  However, if the
  functions $b, \sigma, f\ \mbox{and}\ \Phi$\ are deterministic it is well known that also $u$ is a
  deterministic function of $(t,x)$\ (see, e.g., Proposition 2.4 in
Peng~\cite{Pe1}). \er From now on, let us suppose that
$$
\begin{array}{ll}
\mbox{(vi)} & \mbox{The coefficients}\ b, \sigma, f\ \mbox{and}\
\Phi\
\mbox{are deterministic, i.e., }\ \mbox{independent of}\\
& (\omega', \omega)\in
\bar{\Omega}=\Omega\times\Omega.\\
\end{array}
\eqno{\mbox{(H5.2)}} $$

 The function $u$\ and the random field $Y^{t,\zeta},\ (t, \zeta)\in [0,
T]\times L^2(\Omega,{\cal{F}}_t,P;{\mathbb{R}}^n),$\ are related by
the following theorem.
 \bt Under the assumptions (H4.1) and (H5.1), for any $t\in [0, T]$\ and $\zeta\in
L^2(\Omega,{\cal{F}}_t,P;{\mathbb{R}}^n),$\ we have \be
u(t,\zeta)=Y_t^{t,\zeta},\ \mbox{ P-a.s.}. \ee \et

Indeed, once $X^{0,x_0}$\ and $Y^{0,x_0}$\ are determined, the
coefficients of SDE (5.1) and BSDE (5.2) are well-determined,
deterministic, depend only on $(t,x)$ and $(t,x,y,z)$, respectively,
and satisfy the standard growth and Lipschitz conditions. Hence,
Theorem 5.1 is a consequence of the corresponding result in
Peng~\cite{Pe1} (Theorem 4.7) (or, see also Theorem 6.1 in Buckdahn
and Li~\cite{BL}). \vskip0.3cm

We now discuss (the generalized) DPP for our FBSDE (5.1), (5.2). For
this end we have to define the family of (backward) semigroups
associated with BSDE (5.2). This notion of stochastic backward
semigroup was first introduced by Peng~\cite{Pe1} and originally
applied to study the DPP for stochastic control problems. Our
approach extends Peng's ideas to the framework of Mean-Field FBSDE.
However, we change the definition of the stochastic backward
semigroup to simplify the proof of the existence of a viscosity
solution of the associated PDE.

 Given the initial data $(t,x)$, a positive number $\delta\leq T-t$\ and a real-valued
 random variable $\eta \in L^2 (\Omega,
{\mathcal{F}}_{t+\delta},P;{\mathbb{R}})$, we put \be G^{t,
x}_{s,t+\delta} [\eta]:= \tilde{Y}_s^{t,x},\ \hskip0.5cm s\in[t,
t+\delta], \ee where the couple $(\tilde{Y}_s^{t,x},
\tilde{Z}_s^{t,x})_{t\leq s \leq t+\delta}$ is the solution of the
following BSDE with the time horizon $t+\delta$: \be \left
\{\begin{array}{rcl}
 -d\tilde{Y}_s^{t,x} \!\!\!& = &\!\!\! E'[f(s,(X_s^{0,x_0})',X_s^{t,x},
(Y_s^{0,x_0})',Y_s^{t,x},\tilde{Z}_s^{t,x})]ds-\tilde{Z}_s^{t,x} dB_s , \hskip 1cm s\in [t,t+\delta],\\
 \tilde{Y}_{t+\delta}^{t,x}\!\!\! & =& \!\!\!\eta;
\end{array}\right.
\ee here $X^{t,x}$\ is the solution of SDE (5.1) and $Y^{t,x}$\ is
the solution of BSDE (5.2). Then, obviously, for the solution
$(Y^{t,x}, Z^{t,x})$\ of BSDE (5.2) we have \be G^{t,x}_{t,T}
[\tilde{\Phi} (X^{t,x}_T)] =G^{t,x}_{t,t+\delta}
[Y^{t,x}_{t+\delta}],\ \ 0\leq t<t+\delta\leq T. \ee

Let us point out that in difference to Peng's definition of the
backward semigroup the driver of our BSDE depends on the processes
$Y^{0,x_0}$ and $Y^{t,x}$ given by BSDE (5.2). The choice to let the
driver $f$\ depend on these processes and not on
$\widetilde{Y}^{0,x_0}$, $\widetilde{Y}^{t,x}$ will simplify the
proof of the existence theorem for the associated nonlocal PDEs in
Section 6.

Moreover, we have the following DPP\be
\begin{array}{rcl}
 u(t,x)& = &Y_t^{t,x}=G^{t,x}_{t,T} [\tilde{\Phi} (X^{t,x}_T)]
  =G^{t,x}_{t,t+\delta} [Y^{t,x}_{t+\delta}]\\
  &=&G^{t,x}_{t,t+\delta} [u(t+\delta,X^{t,x}_{t+\delta})],
\end{array}
\ee whose simple form explains by the fact that our stochastic
evolution system doesn't depend on  a control. Here, for the latter
relation we have used, that due to the uniqueness of the solution of
SDE (5.1) and of BSDE (5.2),
$Y^{t,x}_{t+\delta}=Y^{t+\delta,X^{t,x}_{t+\delta}}_{t+\delta},\ \
\mbox{P-a.s.},$\ so that then from Theorem 5.1 it follows that
$Y^{t,x}_{t+\delta}=u(t+\delta,X^{t,x}_{t+\delta}),\ \
\mbox{P-a.s.}$
 \br If $f$\ is independent of $(y, z)$\ it holds that
$$G^{t,x}_{s,t+\delta}[\eta]=E[\eta + \int_s^{t+\delta}
E'[f(r,(X^{0,x_0}_r)',X^{t,x}_r, (Y^{0,x_0}_r)')]dr|{\cal{F}}_s],\ \
s\in [t, t+\delta],\ \mbox{P-a.s.}$$ \er

In (5.5) we have already seen that the value function $u(t,x)$\ is
Lipschitz continuous in $x$, uniformly in $t$. Relation (5.10) now
allows also to study the continuity property of $u(t,x)$\ in $t$.
 \bt\mbox{ }Let us suppose that the assumptions (H4.1), (H5.1) and (H5.2)
hold. Then the value function $u(t,x)$ is
 $\frac{1}{2}-$H\"{o}lder continuous in $t$, locally uniformly with respect to $x$: There exists a constant C such that,
  for every $x\in {\mathbb{R}}^n,\ t,\ t'\in [0, T]$,
  $$
|u(t, x)-u(t', x)|\leq C(1+|x|)|t-t'|^{\frac{1}{2}}.
  $$
  \et
\noindent \textbf{Proof}.  Let $(t, x)\in [0,T]\times
{\mathbb{R}}^n$\ and $\delta>0$\ be arbitrarily given such that
$0<\delta\leq T-t$. Our objective is to prove the following
inequality by using DPP:
 \be -C(1+|x|)\delta^{\frac{1}{2}}\leq u(t,x)-u(t+\delta ,x)\leq
C(1+|x|)\delta^{\frac{1}{2}}.
 \ee
 From it we obtain immediately that $u$ is $\frac{1}{2}-$H\"{o}lder continuous in
 $t$. We will only check the second inequality in (5.11), the
 first one can be shown in a similar way. To this end we note that
 due to (5.10),
\be u(t,x)-u(t+\delta ,x) = I^1_\delta +I^2_\delta, \ee
 where
$$
\begin{array}{lll}
I^1_\delta & := & G^{t,x}_{t,t+\delta}[u(t+\delta,
X^{t,x}_{t+\delta})]
                   -G^{t,x}_{t,t+\delta} [u(t+\delta,x)], \\
I^2_\delta & := & G^{t,x}_{t,t+\delta} [u(t+\delta,x)] -u(t+\delta
,x).
\end{array}
$$
From Lemmata 2.3 and 3.1 and the estimate (5.5) we obtain that,
$$
\begin{array}{rcl}
|I^1_\delta | &\leq& [CE(|u(t+\delta ,X^{t,x}_{t+\delta})
                  -u(t+\delta ,x)|^2|{{\mathcal{F}}_t})]^{\frac{1}{2}}\\
              & \leq&[CE(|X^{t,x}_{t+\delta} -x|^2|{{\mathcal{F}}_t})]^{\frac{1}{2}},
\end{array}
$$
and since $E[|X^{t,x}_{t+\delta} -x|^2|{{\mathcal{F}}_t}] \leq
C(1+|x|^2) \delta $ we deduce that $|I^1_\delta| \leq C
(1+|x|)\delta^{\frac{1}{2}}$. From the definition of
$G^{t,x}_{t,t+\delta}[\cdot]$\ (see (5.7)) we know that the second
term $I^2_\delta $ can be written as£º
$$
\begin{array}{llll}
I^2_\delta  & = & E[u(t+\delta ,x) +\int^{t+\delta}_t
E'[f(s,(X_s^{0,x_0})',X_s^{t,x}, (Y_s^{0,x_0})', Y_s^{t,x},
\tilde{Z}_s^{t,x})]ds  \\
 &         &   -\int^{t+\delta}_t \tilde{Z}^{t,x}_s dB_s|{{\mathcal{F}}_t}] -u(t+\delta ,x)  \\
 &    =  &  E[\int^{t+\delta}_t E'[f(s,(X_s^{0,x_0})',X_s^{t,x}, (Y_s^{0,x_0})', Y_s^{t,x},
\tilde{Z}_s^{t,x})]ds|{{\mathcal{F}}_t}].
\end{array}
$$
Then, with the help of the Schwartz inequality, and the estimates
(4.2), (5.3)-(i) for the BSDEs (5.2) and (5.8) (with
$\eta=u(t+\delta,x)$)\ and (5.5) we have
$$
\begin{array}{lll}
|I^2_\delta | & \leq \delta^{\frac{1}{2}}
     E[\int^{t+\delta}_t |E'[f(s,(X_s^{0,x_0})',X_s^{t,x}, (Y_s^{0,x_0})', Y_s^{t,x},
\tilde{Z}_s^{t,x})]|^2ds|{{\mathcal{F}}_t}]^{\frac{1}{2}}  \\
& \leq\delta^{\frac{1}{2}}E[\int^{t+\delta}_t
(|E'[f(s,(X_s^{0,x_0})',X_s^{t,x}, (Y_s^{0,x_0})', 0, 0)]|+C|
Y^{t,x}_s|+C|\tilde{Z}^{t,x}_s|)^2ds|{{\mathcal{F}}_t}]^{\frac{1}{2}}\\
& \leq C\delta^{\frac{1}{2}}E[\int^{t+\delta}_t (|1+|X^{t,x}_s|+|
Y^{t,x}_s|
+|\tilde{Z}^{t,x}_s|)^2ds|{{\mathcal{F}}_t}]^{\frac{1}{2}}\\
 & \leq C (1+|x|)\delta^{\frac{1}{2}}.
\end{array}
$$
Hence, from (5.12), we get the second inequality of
(5.11)
$$u(t,x)-u(t+\delta ,x) \leq C (1+|x|)\delta^{\frac{1}{2}}.$$ The proof is complete.\endpf

\section{\large Viscosity Solution Of PDE: Existence Theorem}

 \hskip1cm In this section we consider the following PDE \be\left
\{\begin{array}{ll}
 &\!\!\!\!\! \frac{\partial }{\partial t} u(t,x) + Au(t,x)
+ E[f(t, X_t^{0, x_0}, x, u(t, X_t^{0, x_0}), u(t,x),
Du(t,x).E[\sigma(t, X_t^{0, x_0},
 x)])]=0,\\
 & \hskip 4.3cm (t,x)\in [0,T)\times {\mathbb{R}}^n ,  \\
 &\!\!\!\!\!  u(T,x) =E[\Phi (X_T^{0, x_0}, x)], \hskip0.5cm   x \in
 {\mathbb{R}}^n,
 \end{array}\right.
\ee with
$$ Au(t,x):=\frac{1}{2}tr(E[\sigma(t, X_t^{0, x_0}, x)]
 E[\sigma(t, X_t^{0, x_0},x)]^{T}D^2u(t,x))+Du(t,x).E[b(t,
 X_t^{0, x_0},x)].$$ Here the functions $b, \sigma, f\ \mbox{and}\ \Phi$\ are
supposed to satisfy (H4.1), (H5.1) and (H5.2), respectively, and
$X^{0, x_0}$\ is the solution of the Mean-Field SDE (5.1).

We attract the reader's attention to the fact that, since
$$\begin{array}{rcl}& & E[f(t, X_t^{0, x_0}, x, u(t, X_t^{0, x_0}), u(t,x), Du(t,x).E[\sigma(t, X_t^{0, x_0},
 x)])]\\
 =& & \int_{{\mathbb{R}}^n}f(t, x', x, u(t, x'), u(t,x), Du(t,x).E[\sigma(t, X_t^{0, x_0},
 x)])P_{X_t^{0, x_0}}(dx'),\end{array}$$
the above equation is in fact a nonlocal PDE.

 In this section we want to prove that
the value function $u(t, x)$ introduced by (5.4) is the
 viscosity solution of equation (6.1). For this we extend Peng's BSDE approach~\cite{Pe1}
developed in the framework of stochastic control theory to that of
the Mean-Field FBSDE. The difficulties related with this extension
come from the fact that now, contrarily to the framework of
stochastic control theory studied by Peng, we have to do with
nonlocal PDEs. Moreover, in difference to~\cite{BBE} the nonlocal
term is not generated by a diffusion process with jumps. This fact
is the source of difficulties mainly in the proof of the uniqueness
of the viscosity solution (given in the next section) which are
different to those in~\cite{BBE}. Let us first recall the definition
of a viscosity solution of equation (6.1). The reader more
interested in viscosity solutions is referred to Crandall, Ishii and
Lions~\cite{CIL}.

\bde\mbox{ } A real-valued
continuous function $u\in C_p([0,T]\times {\mathbb{R}}^n )$ is called \\
  {\rm(i)} a viscosity subsolution of equation (6.1) if, firstly, $u(T,x) \leq E[\Phi (X_T^{0, x_0}, x)],\ \mbox{for all}\ x \in
  {\mathbb{R}}^n$, and if, secondly, for all functions $\varphi \in C^3_{l, b}([0,T]\times
  {\mathbb{R}}^n)$ and $(t,x) \in [0,T) \times {\mathbb{R}}^n$ such that $u-\varphi $\ attains its
  local maximum at $(t, x)$,
     $$\begin{array}{ll}
     &\frac{\partial }{\partial t} \varphi(t,x) + D\varphi(t,x).E[b(t, X_t^{0, x_0}, x)]+\frac{1}{2}tr(E[\sigma(t, X_t^{0, x_0}, x)]
 E[\sigma(t, X_t^{0, x_0},x)]^{T}D^2\varphi(t,x))\\
 &+ E[f(t, X_t^{0, x_0}, x, u(t, X_t^{0, x_0}), u(t,x), D\varphi(t,x).E[\sigma(t, X_t^{0, x_0},
 x)])]\geq 0;\end{array}
     $$
{\rm(ii)} a viscosity supersolution of equation (6.1) if, firstly,
$u(T,x) \geq E[\Phi (X_T^{0, x_0}, x)],\ \mbox{for all}\ x \in
  {\mathbb{R}}^n$, and if, secondly, for all functions $\varphi \in C^3_{l, b}([0,T]\times
  {\mathbb{R}}^n)$ and $(t,x) \in [0,T) \times {\mathbb{R}}^n$ such that $u-\varphi $\ attains its
  local minimum at $(t, x)$,
     $$\begin{array}{ll}
     &\frac{\partial }{\partial t} \varphi(t,x) + D\varphi(t,x).E[b(t, X_t^{0, x_0}, x)]+\frac{1}{2}tr(E[\sigma(t, X_t^{0, x_0}, x)]
 E[\sigma(t, X_t^{0, x_0},x)]^{T}D^2\varphi(t,x))\\
 &+ E[f(t, X_t^{0, x_0}, x, u(t, X_t^{0, x_0}), u(t,x), D\varphi(t,x).E[\sigma(t, X_t^{0, x_0},
 x)])] \leq 0;\end{array}
     $$
 {\rm(iii)} a viscosity solution of equation (6.1) if it is both a viscosity sub- and a supersolution of equation
     (6.1).\ede
\br {\rm(i)} $C_p([0,T]\times R^n)=\{u\in C([0,T]\times R^n):
\mbox{There exists some constant}\ p>0\ \mbox{such that}$ $
\sup_{(t,x)\in[0, T]\times R^n}\frac{|u(t,x)|}{1+|x|^p}<
+\infty\}.$\\
{\rm(ii)} Usually the definition of a viscosity solution is given
for test fields $\varphi$\ of class $C^{1,2}$. However, it can be
shown that it is sufficient to work with test functions from
$C^3_{l, b}([0,T]\times {\mathbb{R}}^n)$. The space $C^3_{l,
b}([0,T]\times {\mathbb{R}}^n)$ denotes the set of the real-valued
functions that are continuously differentiable up to the third order
and whose derivatives of order from 1 to 3 are bounded.\er

We now can state the main statement of this section.
 \bt Under the assumptions (H4.1), (H5.1) and (H5.2) the function $u(t,x)$\ defined by (5.4) is a viscosity solution of
equation (6.1). \et

The proof of the theorem uses the BSDE method of Peng~\cite{Pe1}.
However, it is simplified by the specific choice of our stochastic
backward semigroup. For the proof of this theorem we need four
auxiliary lemmata. To abbreviate notations we put, for some
arbitrarily chosen but fixed $\varphi \in C^3_{l, b} ([0,T] \times
{\mathbb{R}}^n)$,
 \be
\begin{array}{lll}
     F(s,x,z)=&\!\!\!\! \frac{\partial }{\partial s}\varphi (s,x) +
     \frac{1}{2}tr(\tilde{\sigma}\tilde{\sigma}^{T}(s,x)D^2\varphi)+ D\varphi.\tilde{b}(s, x) \\
        &+ \tilde{f}^{u}(s, x, z+ D\varphi (s,x).\tilde{\sigma}(s,x)), \\
     \end{array}
\ee where
$$\begin{array}{rcl}\tilde{\sigma}(s,x)&=& E[\sigma(s,
X_s^{0, x_0}, x)],\ \tilde{b}(s,x)=E[b(s, X_s^{0, x_0}, x)];\\
\tilde{f}^{u}(s, x ,z)&=& E[f(s,X_s^{0, x_0}, x, u(s, X_s^{0,
x_0}),u(s,x), z)],\end{array}$$ $(s,x, z)\in [0,T] \times
{\mathbb{R}}^n \times {\mathbb{R}}^d,$\ and we consider the
following BSDE defined on the interval $[t,t+\delta]\ (0<\delta\leq
T-t):$
    \be \left \{\begin{array}{rl}
      -dY^{1}_s =&\!\!\!\! F(s,X^{t,x}_s,Z^{1}_s)ds
                   -Z^{1}_s dB_s, \\
     Y^{1}_{t+\delta}=&\!\!\!\!0,
     \end{array}\right.
     \ee
     where the process $X^{t,x}$\ has been introduced by equation
     $(5.1)$.
\br\mbox{}It's not hard to check that $F(s,X^{t,x}_s,z)$\ satisfies
(A1) and (A2). Thus, due to Lemma 2.1 equation (6.3) has a unique
solution. \er

      We can characterize the solution process $Y^{1}$ as follows:
    \bl\mbox{  } For every $s\in [t,t+\delta]$, we have the following
relationship:
    \be
     Y^{1}_s = G^{t,x}_{s,t+\delta} [\varphi (t+\delta ,X^{t,x}_{t+\delta})]
                -\varphi (s,X^{t,x}_s), \hskip 0.5cm
              \mbox{P-a.s.}\ee
 \el
\noindent \textbf{Proof}. We recall that $G^{t,x}_{s,t+\delta}
[\varphi (t+\delta, X^{t,x}_{t+\delta})]$ is defined with the help
of the solution of the BSDE
     $$
     \left \{\begin{array}{rl}
     -d\tilde{Y}_s =\!\!\! & E'[f(s,(X^{0,x_0}_s)',X^{t,x}_s, (Y^{0,x_0}_s)', Y_s^{t,x},\tilde{Z}_s)]ds
      -\tilde{Z}_s dB_s , \hskip 0.2cm s\in [t,t+\delta], \\
     \tilde{Y}_{t+\delta}=\!\!\!& \varphi (t+\delta
     ,X^{t,x}_{t+\delta}),
     \end{array}\right.
     $$
by the following formula:
     \be
     G^{t,x}_{s,t+\delta} [\varphi (t+\delta ,X^{t,x}_{t+\delta})] =\tilde{Y}_s, \hskip 0.5cm
     s\in [t,t+\delta]  \ee
(see (5.7)). We also recall that $E'[f(s,(X_s^{0, x_0})', X^{t,x}_s,
(Y_s^{0, x_0})', Y_s^{t,x},\tilde{Z}_s)]=\tilde{f}^{u}(s, X^{t,x}_s
,\tilde{Z}_s)$. Therefore, we only need to prove that
$\tilde{Y}_s-\varphi (s,X^{t,x}_s)\equiv Y^{1}_s.$\ This result can
be obtained easily by applying It$\hat{o}$'s formula to $\varphi
(s,X^{t,x}_s)$. Indeed, we get that the stochastic differentials of
$\tilde{Y}_s -\varphi (s,X^{t,x}_s)$\ and $Y^{1}_s$\ coincide, while
at the terminal time $t+\delta$, $\tilde{Y}_{t+\delta} - \varphi
(t+\delta ,X^{t,x}_{t+\delta}) =0 = Y^{1}_{t+\delta}.$\ So the proof
is complete.\endpf \vskip 0.3cm

We now introduce the deterministic function
$$Y_s^2=\int_s^{t+\delta}F(r,x,0)dr,\ s\in[t, t+\delta].$$
Obviously, the couple $(Y^2, Z^2)=(Y^2, 0)$\ is the unique solution
of the following (deterministic) BSDE in which the driving process
$X^{t,x}$ is replaced by its deterministic initial value $x$: \be
     \left \{\begin{array}{rl}
     -dY^{2}_s=&\!\!\!  F(s,x, Z^{2}_s)ds -Z^{2}_s dB_s,
         \\
     Y^{2}_{t+\delta}=&\!\!\! 0,
         \hskip 0.3cm s\in [t,t+\delta].
     \end{array}\right.
 \ee
The following lemma will allow us to neglect the difference
$|Y^{1}_t-Y^{2}_t|$ for sufficiently small $\delta >0$.

\bl We have
 \be |Y^{1}_t-Y^{2}_t| \leq C\delta^{\frac{3}{2}},\ \ \mbox{P-a.s.}
       \ee
\el
 \noindent \textbf{Proof}. We recall that from (4.2) with
 $\zeta=x$\ it follows that there is some constant $C_p$\ depending on $p,\ x,$\
 but not on $\delta>0$, such that
      \be
       E [\sup \limits_{t\leq s \leq t+\delta} |X^{t,x}_s -x|^p|{{\mathcal{F}}_t}] \leq
       C_p\delta^{\frac{p}{2}}.
      \ee
   We now apply Lemma 2.3 combined with (6.8) to the equations (6.3) and (6.6). For this we
set in Lemma 2.3:
       $$\xi_1 =\xi_2 =0,\ g(s,y,z) =g(s,z) =F(s,X^{t,x}_s, z),$$
       $$\varphi_1(s)=0,\ \varphi_2(s)=F(s,x,Z_s^{2})-F(s,X^{t,x}_s,Z_s^{2}).
       $$
       Obviously, the function $g$\ is Lipschitz with respect to $z$, and
       $|\varphi_2(s)|\leq C(1+|x|^2)(|X^{t,x}_s -x|+|X^{t,x}_s -x|^3),$\
       for $s\in [t, t+\delta], (t, x)\in [0, T)\times
       {\mathbb{R}}^n$.
       Thus, with the notation $\rho_0 (r) =(1+|x|^2)(r+r^3),\ r\geq 0, $\ we have
      $$
       \begin{array}{lll}
       &  & E[\int^{t+\delta}_t |Z^{1}_s -Z^{2}_s|^2ds|{\mathcal{F}}_{t}]  \\
       & \leq & CE[\int^{t+\delta}_t \rho^2_0 (|X^{t,x}_s -x|) ds|{\mathcal{F}}_{t}]\\
       & \leq &C \delta E[\sup \limits_{t\leq s \leq t+\delta}\rho^2_0 (|X^{t,x}_s
       -x|)|{\mathcal{F}}_{t}]\\
       &\leq & C\delta^2.
       \end{array}
     $$
Therefore,
       $$
       \begin{array}{llll}
      && |Y^{1}_t -Y^{2}_t|  =|E[(Y^{1}_t -Y^{2}_t )|{\mathcal{F}}_{t}]|  \\
       & = & |E[\int^{t+\delta}_t (F(s,X^{t,x}_s,Z^{1}_s)
                -F(s,x,Z^{2}_s)) ds|{\mathcal{F}}_{t}]|   \\
       & \leq & CE [\int^{t+\delta}_t [\rho_0 (|X^{t,x}_s -x|)+|Z^{1}_s -Z^{2}_s|] ds|{\mathcal{F}}_{t}]  \\
       & \leq & CE[\int^{t+\delta}_t\rho_0 (|X^{t,x}_s -x|)ds|{\mathcal{F}}_{t}] + C\delta^{\frac{1}{2}}
       E[\int^{t+\delta}_t|Z^{1}_s -Z^{2}_s|^2ds|{\mathcal{F}}_{t}]^{\frac{1}{2}} \\
       & \leq & C\delta^{\frac{3}{2}}.
       \end{array}
       $$
   Thus, the proof is
   complete.\endpf
     \vskip 0.3cm

 Now we are able to give the proof of Theorem 6.1:

\noindent \textbf{Proof}. Obviously, $u(T,x)=E[\Phi ( X_T^{0,x_0},
x)],\ x\in {\mathbb{R}}^n $. Let us show that $u$ is a viscosity
supersolution (respectively, subsolution). For this we suppose that
$\varphi \in C^3_{l,b} ([0,T] \times {\mathbb{R}}^n)$\ and $(t,x)\in
[0, T)\times {\mathbb{R}}^n$\ are such that $u-\varphi$\ attains its
minimum (respectively, maximum) at $(t,x).$\ Notice that we can
replace the condition of a local minimum (respectively, maximum) by
that of a global one in the definition of the viscosity
supersolution (respectively, subsolution) since $u$ is continuous
and of at most linear growth. Moreover, without loss of generality
we may also suppose that $\varphi (t,x)=u(t,x)$.\ Then, due to the
DPP (see (5.10)),
     $$
     \varphi (t,x) =u(t,x) =G^{t,x}_{t,t+\delta}
[u(t+\delta, X^{t,x}_{t+\delta})],\ 0\leq\delta\leq T-t,
     $$
 and from $u\geq \varphi$\ (respectively, $u\leq \varphi$) and the
 monotonicity property of $G^{t,x}_{t,t+\delta}[\cdot]$\
 (see Lemma 2.2 and Theorem 3.2)  we obtain
     $$
    G^{t,x}_{t,t+\delta}
[\varphi(t+\delta, X^{t,x}_{t+\delta})] -\varphi (t,x)\leq 0 \
(\mbox{respectively},\  \geq 0),\ \  \mbox{P-a.s.}
     $$
 Thus, from Lemma 6.1,
    $$
       Y^{1}_t \leq 0 \ (\mbox{respectively},\  \geq 0),\ \
\mbox{P-a.s.},
    $$
and furthermore, from Lemma 6.2 we have
     $$
     \int_t^{t+\delta}F(s,x,0)ds=Y^{2}_t \leq
C\delta^{\frac{3}{2}} \ (\mbox{respectively},\  \geq
-C\delta^{\frac{3}{2}}),\ \ \mbox{P-a.s.}
     $$
It then follows that
     $$
      F(t,x, 0)\leq 0 \ (\mbox{respectively},\  \geq 0)
     $$
and from the definition of $F$\ we see that $u$ is a viscosity
 supersolution (respectively, subsolution) of equation (6.1). Finally, we prove that $u$ is a viscosity
 solution of equation (6.1).\endpf

\section{\large Viscosity Solution of PDE: Uniqueness Theorem }
\hskip1cm The objective of this section is to study the uniqueness
of the viscosity solution of PDE (6.1), \be\left \{\begin{array}{ll}
 &\!\!\!\!\! \frac{\partial }{\partial t} u(t,x) + Au(t,x)
+ E[f(t, X_t^{0, x_0}, x, u(t, X_t^{0, x_0}), u(t,x),
Du(t,x).E[\sigma(t, X_t^{0, x_0},
 x)])]=0,\\
 & \hskip 4.3cm (t,x)\in [0,T)\times {\mathbb{R}}^n ,  \\
 &\!\!\!\!\!  u(T,x) =E[\Phi (X_T^{0, x_0}, x)], \hskip0.5cm   x \in
 {\mathbb{R}}^n,
 \end{array}\right.
\ee with
$$ Au(t,x):=\frac{1}{2}tr(E[\sigma(t, X_t^{0, x_0}, x)]
 E[\sigma(t, X_t^{0, x_0},x)]^{T}D^2u(t,x))+Du(t,x).E[b(t,
 X_t^{0, x_0},x)].$$  Here the functions $b, \sigma, f\ \mbox{and}\ \Phi$\ are
supposed to satisfy (H4.1), (H5.1) and (H5.2), respectively.

 We will prove the uniqueness for equation (7.1) in the space $C_p([0, T]\times {\mathbb{R}}^n)$ of
 continuous functions with at most polynomial growth. In an earlier work Barles, Buckdahn
 and Pardoux~\cite{BBE} introduced the space
 of continuous functions

 $\Theta=\{ \varphi\in C([0, T]\times {\mathbb{R}}^n): \exists\ \widetilde{A}>0\
 \mbox{such
 that}\ \lim_{|x|\rightarrow \infty}|\varphi(t, x)|
 \exp\{-\widetilde{A}[\log((|x|^2+1)^{\frac{1}{2}})]^2\}=0,\
 \mbox{uniformly in}\ t\in [0, T]\}.$  \vskip 0.1cm
\noindent Its growth condition is slightly weaker than the
assumption of polynomial growth but more restrictive than that of
exponential growth. They proved in $\Theta$\ the uniqueness of the
viscosity solution of an integro-partial differential equation
associated with a decoupled FBSDE with jumps. It was shown
in~\cite{BBE} that this kind of growth condition is optimal for the
uniqueness. However, as the following example shows, we cannot hope
to have this class $\Theta$\ for the uniqueness also for our type of
PDE.

\noindent{\bf Example 7.1}\ Let $n=d=1,\ \sigma(s,x',x)=\sigma x,\
b(s,x',x)=\frac{\sigma^2}{2} x$\ and $x_0=1$. Then
$X_s^{0,x_0}=\exp\{\sigma B_s\},\ s\in[0,T],$\ and for
$f(s,x',x,y',y,z)=y'$, PDE (7.1) takes the form $$\left
\{\begin{array}{ll}
 &\!\!\!\!\! \frac{\partial }{\partial t} u(t,x) + \frac{\sigma^2}{2}x\frac{\partial }{\partial x}u(t,x)
 +\frac{\sigma^2}{2}x^2\frac{\partial^2 }{\partial x^2}u(t,x)+ E[u(t, \exp\{\sigma B_t\})]=0, \hskip 0.5cm   (t,x)\in [0,T)\times {\mathbb{R}},  \\
 &\!\!\!\!\!  u(T,x) =E[\Phi (\exp\{\sigma B_T\}, x)], \hskip0.5cm   x \in
 {\mathbb{R}}.
 \end{array}\right.
$$
But, for the growth function
$$\tilde{\chi}(t,x)=\exp\{\widetilde{A}[\log((|x|^2+1)^{\frac{1}{2}})]^2\},\ (t,x)\in [0,T)\times
{\mathbb{R}},
$$
which comes from the definition of $\Theta$,\ we have
$$\begin{array}{ll}
 &E[\tilde{\chi}(t,X_t^{0,x_0})]=E[\exp\{\widetilde{A}[\log((|\exp\{\sigma
 B_t\}|^2+1)^{\frac{1}{2}})]^2\}]\\
 &\geq E[\exp\{\widetilde{A}\sigma^2B_t^2\}]\\
 &=+\infty, \ \ \mbox{if}\ t\in [\frac{1}{2\widetilde{A}\sigma^2}, T].\end{array}
$$
As the example shows, for a function $u\in \Theta,$\ the coefficient
$E[f(t, X_t^{0, x_0}, x, u(t, X_t^{0, x_0}), y, z)]$\ may be not
well defined. This is the reason why we restrict the study of the
uniqueness to the smaller class $C_p([0, T]\times {\mathbb{R}}^n)$.

\bt\mbox{ } We assume that (H4.1), (H5.1) and (H5.2) hold. Let $u_1$
(resp., $u_2$) $\in C_p([0, T]\times {\mathbb{R}}^n)$\ be a
viscosity subsolution (resp., supersolution) of equation (7.1). Then
we have
    $$
     u_1 (t,x) \leq u_2 (t,x) , \hskip 0.5cm \mbox{for all}\ \ (t,x) \in [0,T] \times {\mathbb{R}}^n .
    $$
\et

The proof of the theorem will be prepared by the following auxiliary
lemmata.

\bl Let $K$\ be a Lipschitz constant of $f(t,x',.,.,.,.)$, uniformly
in $(t, x')$\ from (H5.1)-{\rm(ii)}, and let $\nu>K$. Then, if $u\in
C_p([0, T]\times {\mathbb{R}}^n)$\ is a viscosity subsolution
(respectively, supersolution) of PDE (7.1) the function
$\bar{u}(t,x)=u(t,x)e^{\nu t},\ (t, x)\in [0, T]\times
{\mathbb{R}}^n,$\ is a viscosity subsolution (respectively,
supersolution) of the following PDE: \be\left\{
 \begin{array}{ll}
 & \frac{\partial }{\partial t} \bar{u}(t,x) + D\bar{u}(t,x).E[b(t, X_t^{0, x_0}, x)]+\frac{1}{2}tr(E[\sigma(t, X_t^{0, x_0}, x)]
 E[\sigma(t, X_t^{0, x_0},x)]^{T}D^2\bar{u}(t,x))\\
 &+ E[\bar{f}(t, X_t^{0, x_0}, x, \bar{u}(t, X_t^{0, x_0}), \bar{u}(t,x),D\bar{u}(t,x).E[\sigma(t, X_t^{0, x_0},
 x)])]=0,\  \hskip 0.5cm   (t,x)\in [0,T)\times {\mathbb{R}}^n;  \\
 & \bar{u}(T,x) =E[\Phi (X_T^{0, x_0}, x)]e^{\nu T}, \hskip0.5cm   x \in {\mathbb{R}}^n,
 \end{array}\right.
\ee where $\bar{f}(t,x',x,y',y,z)=e^{\nu t}f(t,x',x,e^{-\nu
t}y',e^{-\nu t}y,e^{-\nu t}z)-\nu y, \ (t,x',x,y',y,z)\in [0,
T]\times {\mathbb{R}}^n\times {\mathbb{R}}^n\times
{\mathbb{R}}\times {\mathbb{R}}\times {\mathbb{R}}^d,$\ conserves
the properties (H5.1) and (H5.2) of $f$\ and is, moreover, strictly
decreasing in $y$:
$$\begin{array}{ll}
&\bar{f}(t,x',x,y',y_1,z)-\bar{f}(t,x',x,y',y_2,z)\leq
-(\nu-K)(y_1-y_2),\\
&\ \hskip6cm \mbox{for all}\ (t,x',x,y',z) \mbox{and all}\ y_1,
y_2\in {\mathbb{R}}\ \mbox{with}\ y_1\geq y_2.
\end{array}
$$
\el The proof of this well-known transformation is straight-forward
and is hence omitted.

\bl\mbox{  } Let $u_1 \in C_p([0, T]\times {\mathbb{R}}^n)$\ be a
viscosity subsolution and $u_2 \in C_p([0, T]\times
{\mathbb{R}}^n)$\ be a viscosity supersolution of equation (7.2).
Then the function $\omega:= u_1-u_2\in C_p([0, T]\times
{\mathbb{R}}^n)$\ is a viscosity subsolution of the equation
    \be
    \left \{
     \begin{array}{lll}
     &\!\!\!\!\!-\nu\omega(t,x)+\frac{\partial }{\partial t} \omega(t,x) + \frac{1}{2}tr(\tilde{\sigma}\tilde{\sigma}^{T}(t, x)D^2\omega)
     + D\omega.\tilde{b}(t, x)+ K|\omega(t,x)| +\\
 &\!\!\!\!\!\mbox{ }\hskip1cm +KE[(\omega(t, X_t^{0,x_0}))^+]+K|D\omega.\tilde{\sigma}(t, x)|= 0, \ \hskip2cm  (t, x)\in [0, T)\times
 {\mathbb{R}}^n,\\
&\!\!\!\!\!\omega(T,x) =0,\ \hskip1cm  x \in {\mathbb{R}}^n.
     \end{array}\right.
   \ee
 (See (6.2) for the notations $\tilde{b}$\ and $\tilde{\sigma}$.) \el
\noindent The proof of this lemma follows from that of Lemma 3.7
in~\cite{BBE}, it turns out to be even simpler because contrary to
Lemma 3.7 in~\cite{BBE} we don't have any integral part generated by
a diffusion process with jumps here in the equations (5.1) and
(5.2).
 \bl\mbox{  }For any
$A>1,$\ there exists $C_1>0$\ such that the function
$$\chi(t,x)=Ae^{C_1(T-t)}\psi(x),$$
with
$$\psi(x)=(|x|^2+1)^{\frac{p}{2}},\ x\in {\mathbb{R}}^n,\ p>1,$$
satisfies
    \be
     \begin{array}{lll}
     &\frac{\partial }{\partial t}\chi(t,x) + \frac{1}{2}tr(\tilde{\sigma}\tilde{\sigma}^{T}(t, x)D^2\chi)+ D\chi.\tilde{b}(t, x)+ K\chi(t,x) +\\
 & K|D\chi(t,x).\tilde{\sigma}(t, x)|+KE[\chi(t,X_t^{0,x_0})]< 0 \ \  \mbox{in}\ [0, T]\times
 {\mathbb{R}}^n.
     \end{array}
    \ee
   \el
\noindent {\bf Proof.} By direct computation we first deduce the
following estimates for the first and second derivatives of $\psi$:
$$ |D\psi(x)|= p\frac{\psi(x)}{|x|^2+1}|x|,\ \ \  |D^2\psi(x)|\leq
p^{2}\frac{\psi(x)}{|x|^2+1},\ \ \ x\in {\mathbb{R}}^n,$$ and
$$\frac{\partial }{\partial t}\chi(t,x)= -C_1\chi(t,x),\ \ \  E[\psi(X_t^{0,x_0})]\leq
C_p\psi(x),\ \ \ x\in {\mathbb{R}}^n.$$ The latter estimate is a
direct consequence of (4.2). Taking into account that the
coefficients $b,\ \sigma$\ are of most linear growth, these
estimates imply for all $t\in [0, T],$
$$  \begin{array}{lll}
     &\frac{\partial }{\partial t}\chi(t,x) + \frac{1}{2}tr(\tilde{\sigma}\tilde{\sigma}^{T}(t, x)D^2\chi)+ D\chi.\tilde{b}(t, x)+ K\chi(t,x) +\\
 &\mbox{}\hskip1cm K|D\chi(t,x).\tilde{\sigma}(t, x)| +KE[\chi(t,X_t^{0,x_0})]\\
 &\leq -\chi(t,x)\{C_1-p^2C-K-KC_p\}\\
 &= -\chi(t,x)< 0,\ \mbox{if}\ C_1:=p^2C+K+KC_p+1 \ \mbox{large
 enough}.
     \end{array}
   $$
\endpf

\noindent Now we can prove the uniqueness theorem.\\
\noindent {\bf Proof of Theorem 7.1.} From Lemma 7.1 we only need to
prove that: If $\bar{u}_1$ (respectively, $\bar{u}_2$) $\in C_p([0,
T]\times {\mathbb{R}}^n)$\ is a viscosity subsolution (respectively,
supersolution) of equation (7.2) then we have
    $$
     \bar{u}_1 (t,x) \leq \bar{u}_2 (t,x) , \hskip 0.5cm \mbox{for all}\ \ (t,x) \in [0,T] \times {\mathbb{R}}^n .
    $$
Let us put $\omega:= \bar{u}_1-\bar{u}_2$. Since $\omega\in C_p([0,
T]\times {\mathbb{R}}^n)$, there exist some constants $C>0$\ and
$q\geq 1$\ such that $|\omega(t,x)|\leq C(1+|x|)^q,\ \ (t, x)\in [0,
T]\times
 {\mathbb{R}}^n.$\ In the definition of $$\chi(t,x)=Ae^{C_1(T-t)}(|x|^2+1)^{\frac{p}{2}}$$
we now choose $p>q$. Then, for any $\alpha
>0$, $\omega(t, x)-\alpha\chi(t, x)$\ is bounded from above in
$[0, T]\times {\mathbb{R}}^n$\ and
$$ M:=\max_{[0, T]\times {\mathbb{R}}^n}(\omega-\alpha\chi)(t, x)e^{-K(T-t)}$$
is achieved at some point $(t_0, x_0)\in [0, T]\times
{\mathbb{R}}^n$\ (depending on $\alpha$). We now have to distinguish between two cases:\\
For the first case we suppose that $\omega(t_0, x_0)\leq 0$, for all $\alpha>0$.\\
Then, obviously $M\leq 0$ and $\bar{u}_1(t, x)-\bar{u}_2(t, x)\leq
\alpha\chi(t, x)$\ in $[0, T]\times {\mathbb{R}}^n$. Consequently,
letting $\alpha$\ tend to zero we obtain
$$\bar{u}_1(t, x)\leq \bar{u}_2(t, x),\ \ \mbox{for all}\ (t, x)\in [0, T]\times {\mathbb{R}}^n. $$
For the second case we assume that there exists some $\alpha>0$\ such that $\omega(t_0, x_0)> 0$.\\
We notice that $\omega(t, x)-\alpha\chi(t,x)\leq (\omega(t_0,
x_0)-\alpha\chi(t_0, x_0))e^{-K(t-t_0)}\ \ \mbox{in}\ \ [0, T]\times
{\mathbb{R}}^n.$\ Then, putting
$$\varphi(t, x)=\alpha\chi(t, x)+(\omega-\alpha\chi)(t_0, x_0)e^{-K(t-t_0)},\ \ (t, x)\in [0, T]\times {\mathbb{R}}^n,$$
we get $\omega-\varphi\leq 0=(\omega-\varphi)(t_0, x_0)\ \mbox{in}\
\ [0, T]\times {\mathbb{R}}^n.$\ Consequently, since $\omega$\ is a
viscosity subsolution of (7.3) from Lemma 7.2,
 $$
     \begin{array}{lll}
     &-\nu\omega(t_0, x_0)+\frac{\partial }{\partial t}\varphi(t_0, x_0) +
     \frac{1}{2}tr(\tilde{\sigma}\tilde{\sigma}^{T}(t_0, x_0)D^2\varphi(t_0, x_0))+ D\varphi(t_0, x_0).\tilde{b}(t_0, x_0)+\\
  & K|\varphi(t_0, x_0)| +
 K|D\varphi(t_0, x_0).\tilde{\sigma}(t_0, x_0)|+KE[(\varphi(t_0, X_{t_0}^{0,x_0}))^+]\geq 0.
     \end{array}
$$
Moreover, due to our assumption that $\omega(t_0, x_0)>0$\ and since
$\omega(t_0, x_0)=\varphi(t_0, x_0)$\ we can replace $K|\varphi(t_0,
x_0)|$\ by $K\varphi(t_0, x_0)$\ in the above formula. Then, from
the definition of $\varphi$\ and Lemma 7.3,
$$
     \begin{array}{lll}
     &0\leq -\nu\omega(t_0, x_0)+\frac{\partial }{\partial t}\varphi(t_0, x_0) +
     \frac{1}{2}tr(\tilde{\sigma}\tilde{\sigma}^{T}(t_0, x_0)D^2\varphi(t_0, x_0))+ D\varphi(t_0, x_0).\tilde{b}(t_0, x_0)+\\
  &\ \hskip0.5cm K|\varphi(t_0, x_0)| +
 K|D\varphi(t_0, x_0).\tilde{\sigma}(t_0, x_0)|+KE[(\varphi(t_0,
 X_{t_0}^{0,x_0}))^+]\\
&=\alpha\frac{\partial }{\partial t}\chi(t_0, x_0) + \alpha K
\chi(t_0, x_0)+
     \alpha\frac{1}{2}tr(\tilde{\sigma}\tilde{\sigma}^{T}(t_0, x_0)D^2\chi(t_0, x_0))+ \alpha D\chi(t_0, x_0).\tilde{b}(t_0, x_0)+\\
  &\ \hskip0.5cm \alpha K|D\chi(t_0, x_0).\tilde{\sigma}(t_0, x_0)|+KE[(\alpha\chi(t_0,
 X_{t_0}^{0,x_0})+(\omega-\alpha\chi)(t_0, x_0))^+]-\nu\omega(t_0, x_0)\\
&\leq  \alpha\{\frac{\partial \chi}{\partial t} (t_0, x_0)
+K\chi(t_0, x_0)+
\frac{1}{2}tr(\tilde{\sigma}\tilde{\sigma}^{T}(t_0, x_0)D^2\chi(t_0, x_0))+ D\chi(t_0, x_0).\tilde{b}(t_0, x_0) \\
 &\ \hskip0.5cm  + K|D\chi(t_0, x_0).\tilde{\sigma}(t_0, x_0)|+KE[\chi(t_0,X_{t_0}^{0,x_0})]\} < 0,
     \end{array}
   $$
since $\nu>K,$\ which is a contradiction. Thus, the proof is
complete.\endpf

\br Obviously, since the value function $u(t,x)$\ is of at most
linear growth it belongs to $C_p([0, T]\times {\mathbb{R}}^n)$, and
so $u(t,x)$ is the unique viscosity solution in $C_p([0, T]\times
{\mathbb{R}}^n)$ of equation (6.1).\er

\end{document}